\documentclass[12pt]{article}
\usepackage{amsmath}
\usepackage{amssymb}
\usepackage{vatola}

\def\q{\quad}
\def\qq{\qquad}
\def\qtq#1{\q\t{#1}\q}
\def\mod#1{\ (\text{\rm mod}\ #1)}
\def\t{\text}
\def\f{\frac}
\def\e{\equiv}
\def\b{\binom}
\def\sls#1#2{(\f{#1}{#2})}
 \def\ls#1#2{\big(\f{#1}{#2}\big)}
\def\Ls#1#2{\Big(\f{#1}{#2}\Big)}
\def\ap{\langle a\rangle_p}

\def\qp#1{q_p(#1)}

\let \pro=\proclaim
\let \endpro=\endproclaim

\begin{document}
 \centerline {\bf
Supercongruences for sums involving $\binom ak^m$}
\par\q\newline
\centerline{Zhi-Hong Sun}\newline
\centerline{School of Mathematics
and Statistics}
\centerline{Huaiyin Normal University}
\centerline{Huaian, Jiangsu 223300, P.R. China} \centerline{Email:
zhsun@hytc.edu.cn} \centerline{Homepage:
http://maths.hytc.edu.cn/szh1.htm}
 \abstract{ \par Let $p$ be an odd prime, and let $a$ be a rational
 $p$-adic integer with $a\not\equiv 0\pmod p$. In this paper,
  using WZ method we establish the congruences for
   $\sum_{k=0}^{p-1} \binom ak^2(-1)^k(1-\frac 2ak)$ modulo $p^2$ and
   $\sum_{k=0}^{p-1} \binom ak^r(1-\frac 2ak)^s$ modulo $p^4$, where
   $r\in\{3,4\}$ and $s\in\{1,3\}$.
 \par\q
 }
 \endabstract
\let\thefootnote\relax \footnotetext {The author is supported by
the National Natural Science Foundation of China (Grant No.
11771173).}

\section*{1. Introduction}
\par The general binomial coefficient $\b ak$ is given by
$$\b a0=1,\q \b a{-k}=0\qtq{and}\b
ak=\f{a(a-1)\cdots(a-k+1)}{k!}\qtq{for}k=1,2,3,\ldots$$ It is well
known that $$\b{-\f 12}k=\f{\b{2k}k}{(-4)^k},\q\b {\f
12}k=-\f{\b{2k}k}{(-4)^k(2k-1)} \qtq{for}k=0,1,2,\ldots$$ and $p\mid
\b{2k}k$ for given odd prime $p$ and $k\in\{\f{p+1}2,\ldots,p-1\}$.
In 1997, Van Hamme [VH] revealed the surprising correspondence
between formulas for some hypergeometric series and congruences for
corresponding partial sums modulo prime powers. Based on known
formulas for certain hypergeometric series involving $\b ak^m$, Van
Hamme posed 10 conjectures on congruences for corresponding partial
sums modulo prime powers. His crucial idea is to replace the Gamma
function $\Gamma(x)$ in the formula for hypergeometric series with
the $p-$adic Gamma function $\Gamma_p(x)$ in the corresponding
congruences, where
$\Gamma(x)=\int_0^{\infty}t^{x-1}\t{e}^{-t}\t{d}t\q (x>0).$ For an
odd prime $p$ and rational $p-$adic integer $x$,
 the $p$-adic Gamma function $\Gamma_p(x)$ is
  defined by
  $$\Gamma_p(0)=1,\q \Gamma_p(n)=(-1)^n\prod_{
  \substack{k\in\{1,2,\ldots,n-1\}\\p\nmid
  k}}k\qtq{for}n=1,2,3,\ldots$$
  and
  $$\Gamma_p(x)=\lim_{\substack{n\in\{0,1,\ldots\}\\ |x-n|_p\rightarrow
  0}}
  \Gamma_p(n).$$
As typical examples, Van Hamme conjectured that for any prime $p>3$,
$$\align &\sum_{k=0}^{p-1}(4k+1)\b{-1/2}k^3\e -\f p{\Gamma_p(\f 12)^2}
\e (-1)^{\f{p-1}2}p\mod {p^3},\tag 1.1
\\&\sum_{k=0}^{p-1}(8k+1)\b{-1/4}k^3\e-\f p{\Gamma_p(\f
14)\Gamma_p(\f 34)}\mod {p^3}\qtq{for}p\e 1\mod 4,\tag 1.2
\\&\sum_{k=0}^{p-1}(8k+1)\b{-1/4}k^4\e \f{p\Gamma_p(\f
12)\Gamma_p(\f 14)}{\Gamma_p(\f 34)}\mod {p^3}\qtq{for}p\e 1\mod
4,\tag 1.3
\endalign$$
which can be viewed as the $p$-analogues of the known formulas:
$$\align &\sum_{k=0}^{\infty}(4k+1)\b{-1/2}k^3=\f 2{\pi}=\f 2{\Gamma(\f
12)^2},
\\&\sum_{k=0}^{\infty}(8k+1)\b{-1/4}k^3=\f{2\sqrt 2}{\pi}=\f
4{\Gamma(\f 14)\Gamma(\f 34)},
\\&\sum_{k=0}^{\infty}(8k+1)\b{-1/4}k^4=\f{2\sqrt 2}{\Gamma(\f
34)^2\sqrt{\pi}}.\endalign$$ We note that (1.1) was first proved by
Mortenson [Mo], and (1.2)-(1.3) were proved by Swisher [Sw] via
hypergeometric series identities and the properties of $p$-adic
Gamma functions.
\par Let $m\ge 2$ be an integer, and let $p>3$ be a prime. Using
Swisher's method, B. He [He1] proved that for $p\e \pm 1\mod m$,
$$\sum_{k=0}^{p-1}(2km+1)\b{-\f 1m}k^3
\e\cases (-1)^{\f{p-1}m}p\mod {p^3}&\t{if $m\mid p-1$,}
\\(-1)^{\f{p-(m-1)}m}(m-1)p\mod {p^3}&\t{if $m\mid p+1$.}
\endcases$$
In a similar way, in [He2] He claimed to show that
$$\sum_{k=0}^{p-1}(2km+1)\b{-\f 1m}k^4\e \cases p\f{\Gamma_p(\f
1m)\Gamma_p(1-\f 2m)}{\Gamma_p(1-\f 1m)}\mod {p^4}&\t{if $m\mid
p-1$,}\\(m-1)p\f{\Gamma_p(\f 1m)\Gamma_p(1-\f 2m)}{\Gamma_p(1-\f
1m)}\mod {p^4}&\t{if $m\mid p+1$.}\endcases$$ However, his formula
for $p\e -1\mod m$ is false. As we will see in Section 5,
$\sum_{k=0}^{p-1}(2km+1)\b{-1/m}k^4\e 0\mod {p^2}$ for $p\e -1\mod
m$.
\par
By [G, (6.8) and (10.1)] and the well known formula
$\Gamma(x)\Gamma(1-x)=\f{\pi}{\sin \pi x}$, we have the following
identities due to Dougall:
$$\align &\sum_{k=0}^{\infty}\Big(1-\f 2ak\Big)\b ak^3=\f{\sin
a\pi}{a\pi}=\f 1{a\Gamma(a)\Gamma(1-a)}\q (a\ge -\f 13),\tag 1.4
\\&\sum_{k=0}^{\infty}\Big(1-\f 2ak\Big)\b ak^4=\f{\sin
a\pi }{a\pi}\cdot\f{\Gamma(2a+1)}{\Gamma(a+1)}=\f{\Gamma(1+2a)}
{a\Gamma(a)\Gamma(1-a)\Gamma(1+a)}\ (a>-\f 12).\tag 1.5\endalign$$
Inspired by (1.4)-(1.5) and the above work due to Van Hamme, Swisher
and He, in this paper we establish general congruences for
$$\align &\sum_{k=0}^{p-1}\Big(1-\f 2ak\Big)\b ak^3,\q
\sum_{k=0}^{p-1}\Big(1-\f 2ak\Big)^3\b ak^3,
\\&\sum_{k=0}^{p-1}\Big(1-\f 2ak\Big)\b ak^4,\q
\sum_{k=0}^{p-1}\Big(1-\f 2ak\Big)^3\b ak^4\endalign$$ modulo $p^4$
via a natural and elementary approach, where $p$ is an odd prime and
$a$ is a rational $p-$adic integer with $a\not\e 0\mod p$. We remark
that Guo [Guo] proved that
$$\sum_{k=0}^{p-1}(4k+1)^3\b{-1/2}k^3\e -3(-1)^{\f{p-1}2}p\mod
{p^3},\tag 1.6$$ and Z.W. Sun [Su1] showed that
$$\sum_{k=0}^{p-1}(4k+1)\b{-1/2}k^3\e
(-1)^{\f{p-1}2}p+p^3E_{p-3}\mod {p^4},\tag 1.7$$ where $\{E_n\}$ are
Euler numbers given by
$$E_0=1,\q E_{2n-1}=0,\q E_{2n}=-\sum_{k=1}^n\b
{2n}{2k}E_{2n-2k}\q(n\ge 1).$$
 As typical results in this paper, we have:
$$\align&\sum_{k=0}^{p-1}\Big(1-\f 2ak\Big)\b ak^3\e
(-1)^{\ap}\f{a-\ap}{a}+\f{(a-\ap)^3}aE_{p-3}(-a)\mod {p^4},\tag 1.8
\\&\sum_{k=0}^{p-1}\Big(1-\f 2ak\Big)^3\b ak^3
\e -3(-1)^{\ap}\f{a-\ap}{a}\tag
1.9\\&\q\qq\qq\qq\qq\qq+\Ls{a-\ap}a^3\big(4-3a^2E_{p-3}(-a)\big)\mod
{p^4},\endalign$$ where $\ap$ is the unique integer in the set
$\{0,1,\ldots,p-1\}$ such that $a\e \ap\mod p$ and $E_n(x)$ is the
Euler polynomial given by
$$ E_n(x)=\f
1{2^n}\sum_{k=0}^n\b nk(2x-1)^{n-k}E_k\ (n\ge 0).$$ Clearly, (1.8)
and (1.9) are vast generalizations of (1.6) and (1.7). For any prime
$p>3$, we also prove that
$$\align &\sum_{k=0}^{p-1}(8k+1)\b{-1/4}k^4
\\&\e\cases (-1)^{\f{p-1}4}\big(2xp-\f{p^2}{2x}-\f{p^3}{8x^3}\big)
\mod {p^4}&\t{if $p=x^2+y^2\e 1\mod 4$ and $4\mid x-1$,}
\\3(-1)^{\f{p+1}4}\f{(2p-1-2^{p-1})p^2}{\binom{(p-1)/2}{(p-3)/4}}
\mod {p^4}&\t{if $p\e 3\mod 4$.}
\endcases\endalign$$
and for $p\e 1,3\mod 8$ and so $p=x^2+2y^2$ $(x,y\in\Bbb Z)$,
$$\align\sum_{k=0}^{p-1}(16k+1)\b{-\f 18}k^4
\e\cases (-1)^{\f y2}p\big(2x-\f{p}{2x}-\f{p^2}{8x^3}\big)\mod
{p^4}&\t{if $8\mid p-1$ and $4\mid
x-1$,}\\3p\big(4y-\f{p}{2y}-\f{p^2}{16y^3}\big)\mod {p^4}&\t{if
$8\mid p-3$ and $4\mid y-1$.}
\endcases\endalign$$

\par Let $p$ be an odd prime and $a\in\Bbb Z_p$ with $a\not\e 0\mod p$.
 In [Su2], Z.W. Sun showed that
$$\sum_{k=0}^{p-1}\b ak^2(-1)^k\e (-1)^{\f{\ap}2}\b a{\f{\ap}2}\mod
{p^2}\qtq{for}2\mid \ap.\tag 1.10$$ In this paper, using WZ method
we show that
$$\sum_{k=0}^{p-1}\b ak^2(-1)^k\e (-1)^{\f{p-\ap}2-1}
\f{2(a-\ap)}{a\b{p-\ap}{(p-\ap)/2}}\mod {p^2}\qtq{for}2\nmid
\ap.\tag 1.11$$ We also establish the congruence for
$\sum_{k=0}^{p-1}\b ak^2(-1)^k(1-\f 2ak)$ modulo $p^2$. In
particular,
$$\sum_{k=0}^{p-1}\b ak^2(-1)^k\Big(1-\f 2ak\Big)\e
(-1)^{\f{\ap}2}2^{2\ap}\f{a-\ap}{a\b{\ap}{\ap/2}}
\mod{p^2}\qtq{for}2\mid \ap.\tag 1.12$$

\par In addition to the above notation,
throughout this paper we use the following notations. For $a\in\Bbb
Z$ and positive odd number $m$ let $\sls am$ denote the Jacobi
symbol. For positive integers $a,b$ and $n$, if $n=ax^2+by^2$ for
some integers $x$ and $y$, we briefly write that $n=ax^2+by^2$. For
a prime $p$ let $\Bbb Z_p$ be the set of rational numbers whose
denominator is not divisible by $p$. For an odd prime and
   $a\in\Bbb Z_p$ set $q_p(a)=(a^{p-1}-1)/p$. Let $[x]$ be the greatest
integer not exceeding $x$. Let $H_0=H_0^{(2)}=0$, $H_n=1+\f
12+\cdots+\f 1n$ and
   $H_n^{(2)}=1+\f 1{2^2}+\cdots+\f 1{n^2}$ $(n\ge 1)$.
The Bernoulli numbers $\{B_n\}$ and Bernoulli polynomials
$\{B_n(x)\}$ are given by
$$B_0=1,\q\sum_{k=0}^{n-1}\b nkB_k=0\
(n\ge 2)\qtq{and}B_n(x)=\sum_{k=0}^n\b nkB_kx^{n-k}\ (n\ge 0).$$
   The sequence
    $\{U_n\}$ is given by
$$U_0=1,\q U_{2n-1}=0,\q  U_{2n}=-2\sum_{k=1}^n\b{2n}{2k}U_{2n-2k}\
(n=1,2,3,\ldots).$$

\section*{2. General congruences for sums involving $\b ak^m$}
\pro{Lemma 2.1} Let $m$ and $n$ be positive integers. Then
$$\align&\sum_{k=0}^{n-1}\b ak^m\sum_{r=0}^{[\f{m-1}2]}\b
m{2r+1}\Big(1-\f 2ak\Big)^{2r+1}=2^{m-1}\b{a-1}{n-1}^m,
\\&\sum_{k=0}^{n-1}(-1)^k\b
ak^m\sum_{r=0}^{[m/2]}\b m{2r}\Big(1-\f
2ak\Big)^{2r}=2^{m-1}(-1)^{n-1}\b{a-1}{n-1}^m .\endalign$$
\endpro
Proof. It is clear that
$$\align &\b{a-1}k^m-\b{a-1}{k-1}^m
\\&=\b ak^m\Big(1-\f ka\Big)^m-\b ak^m\Ls ka^m=\b
ak^m\Big(\Ls{1+(1-\f 2ak)}{2}^m-\Ls{1-(1-\f 2ak)}{2}^m\Big)
\\&=\f 1{2^m}\b ak^m\sum_{s=0}^m\b ms(1-(-1)^s)
\Big(1-\f 2ak\Big)^s\\&=\f 1{2^{m-1}}\b
ak^m\sum_{r=0}^{[\f{m-1}2]}\b m{2r+1}\Big(1-\f
2ak\Big)^{2r+1}.\endalign$$ Thus,
$$\align &\sum_{k=0}^{n-1}\b
ak^m\sum_{r=0}^{[\f{m-1}2]}\b m{2r+1}\Big(1-\f 2ak\Big)^{2r+1}
\\&=2^{m-1}\sum_{k=0}^{n-1}\Big(\b{a-1}k^m-\b{a-1}{k-1}^m\Big)
=2^{m-1}\b{a-1}{n-1}^m.\endalign$$ Similarly,
$$\align &(-1)^k\b{a-1}k^m-(-1)^{k-1}\b{a-1}{k-1}^m
\\&=(-1)^k\b ak^m\Big(\Big(1-\f ka\Big)^m+\b ak^m\Big)\\&=(-1)^k\b
ak^m\Big(\Ls{1+(1-\f 2ak)}{2}^m+\Ls{1-(1-\f 2ak)}{2}^m\Big)
\\&=\f 1{2^m}(-1)^k\b ak^m\sum_{s=0}^m\b ms(1+(-1)^s)
\Big(1-\f 2ak\Big)^s\\&=\f 1{2^{m-1}}(-1)^k\b
ak^m\sum_{r=0}^{[m/2]}\b m{2r}\Big(1-\f 2ak\Big)^{2r}.\endalign$$
Thus,
$$\align &\sum_{k=0}^{n-1}(-1)^k\b
ak^m\sum_{r=0}^{[m/2]}\b m{2r}\Big(1-\f 2ak\Big)^{2r}
\\&=2^{m-1}\sum_{k=0}^{n-1}\Big((-1)^k\b{a-1}k^m-(-1)^{k-1}\b{a-1}{k-1}^m\Big)
=2^{m-1}(-1)^{n-1}\b{a-1}{n-1}^m.\endalign$$ This completes the
proof.

\pro{Theorem 2.1} Let $m$ be a positive integer, and let $p$ be an
odd prime. Assume that $a\in\Bbb Z_p$ and $a\not\e 0\mod p$. Then
$$\align&\sum_{k=0}^{p-1}\b ak^m\sum_{r=0}^{[\f{m-1}2]}\b
m{2r+1}\Big(1-\f 2ak\Big)^{2r+1}\e \sum_{k=0}^{p-1}(-1)^k\b
ak^m\sum_{r=0}^{[m/2]}\b m{2r}\Big(1-\f 2ak\Big)^{2r}
\\&\e
2^{m-1}\f{(a-\ap)^m}{a^m}\Big(1+mpH_{\ap}\Big)
\mod{p^{m+2}}.\endalign$$ In particular, taking $m=1,2,3,4$ we have
$$\align &\sum_{k=0}^{p-1}\b ak\Big(1-\f
2ak\Big)\e\sum_{k=0}^{p-1}(-1)^k\b ak\e
\f{a-\ap}a\big(1+pH_{\ap}\big)\mod {p^3},
\\& 2\sum_{k=0}^{p-1}\b ak^2\Big(1-\f
2ak\Big)\e \sum_{k=0}^{p-1}(-1)^k\b ak^2\Big(1+\Big(1-\f
2ak\Big)^2\Big)\\&\q\qq\qq\qq\q\ \;\e
2\f{(a-\ap)^2}{a^2}\big(1+2pH_{\ap}\big)\mod {p^4},
\\&\sum_{k=0}^{p-1}\b ak^3\Big(3\Big(1-\f 2ak\Big)+\Big(1-\f
2ak\Big)^3\Big)\e \sum_{k=0}^{p-1}(-1)^k\b ak^3\Big(1+3\Big(1-\f
2ak\Big)^2\Big)
\\&\q\qq\qq\qq\qq\qq\qq\qq\ \,\e 4\f{(a-\ap)^3}{a^3}\big(1+3pH_{\ap}\big)\mod {p^5}.
\\&4\sum_{k=0}^{p-1}\b ak^4\Big(\Big(1-\f 2ak\Big)+\Big(1-\f
2ak\Big)^3\Big)\e \sum_{k=0}^{p-1}(-1)^k\b ak^4\Big(1+6\Big(1-\f
2ak\Big)^2+\Big(1-\f 2ak\Big)^4\Big)
\\&\q\qq\qq\qq\qq\qq\qq\qq\ \,
\e 8\f{(a-\ap)^4}{a^4}\big(1+4pH_{\ap}\big)\mod {p^6}.
\endalign$$
\endpro
Proof. Set $t=(a-\ap)/p$. By Lemma 2.1,
$$\align&\sum_{k=0}^{p-1}\b ak^m\sum_{r=0}^{[\f{m-1}2]}\b
m{2r+1}\Big(1-\f 2ak\Big)^{2r+1}=\sum_{k=0}^{p-1}(-1)^k\b
ak^m\sum_{r=0}^{[m/2]}\b m{2r}\Big(1-\f 2ak\Big)^{2r}
\\&=2^{m-1}\b{a-1}{p-1}^m=2^{m-1}\b{\ap+pt-1}{p-1}^m.\endalign$$
By [S5, Lemma 4.2],
$$\align \b{\ap+pt-1}{p-1}&\e
\f{pt}{\ap}-\f{p^2t^2}{\ap^2}+\f{p^2t}{\ap}H_{\ap} \e
\f{pt}{a-pt}-\f{p^2t^2}{a^2}+\f{p^2t}aH_{\ap}
\\&\e
\f{pt(a+pt)}{a^2}-\f{p^2t^2}{a^2}+\f{p^2t}aH_{\ap}=\f{pt}a(1+pH_{\ap})
\mod {p^3}.\endalign$$ Hence,
$$\f 1p\b{a-1}{p-1}\e \f{a-\ap}{ap}\big(1+pH_{\ap}\big)\mod {p^2}.\tag 2.1$$
This yields
$$\b{a-1}{p-1}^m\e
p^m\Ls{a-\ap}{ap}^m(1+pH_{\ap})^m\e \Ls{a-\ap}{a}^m(1+mpH_{\ap})
\mod {p^{m+2}}.$$
 Now, putting all the above together yields the result.
 \par{\bf Remark 2.1} Suppose that $p$ is an odd prime. It is well
 known (see for example [S7]) that $H_k\e H_{p-1-k}\mod p$ for
 $k=0,1,\ldots,p-1$ and
 $$\align &H_{\f{p-1}2}\e -2\qp 2\mod p,\q H_{[\f p4]}\e -3\qp 2\mod p,
\\&H_{[\f p3]}\e -\f 32q_p(3)\mod p,\q H_{[\f p6]} \e -2q_p(2)-\f
32q_p(3)\mod p\q (p>3).
\endalign$$

\pro{Corollary 2.1} Let $p$ be an odd prime. Then
$$\align
&\sum_{k=0}^{p-1}(4k+1)\f{\b{2k}k}{(-4)^k}\e (3-2^p)p\mod {p^3},
\\&\sum_{k=0}^{p-1}(4k+1)\f{\b{2k}k^2}{16^k}\e
\sum_{k=0}^{p-1}(8k^2+4k+1)\f{\b{2k}k^2}{(-16)^k} \e
(5-2^{p+1})p^2\mod {p^4},
\\&\sum_{k=0}^{p-1}\big(3(4k+1)+(4k+1)^3\big)\f{\b{2k}k^3}{(-64)^k}
\e
\sum_{k=0}^{p-1}\big(1+3(4k+1)^2\big)\f{\b{2k}k^3}{64^k}\\&\qq\qq\qq\qq\qq\qq\qq\
\;\,\e 4\big(7-6\cdot 2^{p-1}\big)p^3\mod {p^5},
\\&4\sum_{k=0}^{p-1}\big((4k+1)+(4k+1)^3\big)\f{\b{2k}k^4}{256^k}
\e
\sum_{k=0}^{p-1}\big(1+6(4k+1)^2+(4k+1)^4\big)\f{\b{2k}k^4}{(-256)^k}\\&\qq\qq
\qq\qq\qq\qq\qq \e 8\big(9-2^{p+2}\big)p^4\mod {p^6}.
 \endalign$$
\endpro
Proof. Since $\b{-1/2}k=\b{2k}k\f 1{4^k}$ and $H_{\f{p-1}2}\e -2\qp
2\mod p$, taking $a=-\f 12$ in Theorem 2.1 yields the result.

\section*{3. Congruences for $\sum_{k=0}^{p-1}\b ak^2(-1)^k$ and
$\sum_{k=0}^{p-1}\b ak^2(-1)^k(1-\f 2ak)$ modulo $p^2$}
\par For $n=1,2,3,\ldots$ let
$$f_n(a)=\sum_{k=0}^{n-1}\b ak^2(-1)^k.$$
Set
$$G(a,k)=\f{(-1)^{k-1}}{(a+1)}(2k^2-(6(a+2)-2)k+(a+2)(5a+7))\b{a+1}{k-1}^2
\ (k\ge 0).$$ Using Maple one can check that
$$(a+2)\b{a+2}k(-1)^k+4(a+1)\b ak^2(-1)^k=G(a,k+1)-G(a,k).$$
Thus,
$$\align &(a+2)f_n(a+2)+4(a+1)f_n(a)\\&
=\sum_{k=0}^{n-1}\Big((a+2)\b{a+2}k(-1)^k+4(a+1)\b
ak^2(-1)^k\Big)=\sum_{k=0}^{n-1}(G(a,k+1)-G(a,k))
\\&=G(a,n)-G(a,0)=G(a,n).\endalign$$
\pro{Lemma 3.1} Suppose that $p$ is an odd prime and $a\in\Bbb Z_p$
with $a+1\not\e 0\mod p$. Then
$$f_p(a)\e \cases -\f{a+2}{4(a+1)}f_p(a+2)\mod {p^2}&\t{if $a+2\not\e 0\mod p$,}
\\\f{p-(a+2)}2\mod{p^2}&\t{if $a+2\e 0\mod p$.}
\endcases$$
\endpro
Proof. For $\ap<p-2$ we see that $\b{a+1}{p-1}\e
\b{\ap+1}{p-1}=0\mod p$ and so $G(a,p)\e 0\mod {p^2}$. Hence,
$(a+2)f_p(a+2)+4(a+1)f_p(a)=G(a,p)\e 0\mod {p^2}$ and so the result
is true. For $\ap=p-2$ we  have
$$\align\b{a+1}{p-1}&=\f{((a+2)-1)((a+2)-2)\cdots((a+2)-(p-1))}{(p-1)!}
\\&\e (-1)^{p-1}(1-(a+2)H_{p-1})\e 1\mod {p^2}\endalign$$ and so
$$\align G(a,p)&=\f 1{a+1}(2p^2-(6(a+2)-2)p+(a+2)(5a+7))\b{a+1}{p-1}^2
\\&\e\f 1{a+1}(2p-3(a+2))\e 3(a+2)-2p\mod {p^2}.\endalign$$
Also,
$f_p(a+2)=1+\sum_{k=1}^{p-1}\f{(a+2)^2}{k^2}\b{a+1}{k-1}^2(-1)^k\e
1\mod {p^2}.$ Thus,
$$f_p(a)=\f{G(a,p)-(a+2)f_p(a+2)}{4(a+1)}\e \f{3(a+2)-2p-(a+2)}{-4}
=\f{p-(a+2)}2\mod{p^2}.$$ This completes the proof.
\par\q
\par Let $p$ be an odd prime, $a\in\Bbb Z_p$ and $2\mid \ap$. Then
clearly
$$\align \b a{\f {\ap}2}&=\b{\ap+(a-\ap)}{\f{\ap}2}
\e\f{\ap(\ap-1)\cdots(\f{\ap}2+1)\big(1+(a-\ap)\sum_{k=\f{\ap}2+1}^{\ap}\f
1k\big)}{\f{\ap}2!}
\\&=\b{\ap}{\ap/2}\big(1+(a-\ap)(H_{\ap}-H_{\f{\ap}2})\big)\mod
{p^2}.\endalign$$ This together with (1.10) gives
$$\sum_{k=0}^{p-1}\b ak^2(-1)^k\e (-1)^{\f{\ap}2}\b{\ap}{\f{\ap}2}
\big(1+(a-\ap)(H_{\ap}-H_{\f{\ap}2})\big)\mod {p^2}\ \t{for}\
2\mid\ap.\tag 3.1$$
\par Now we turn to present the result in the case $\ap\e 1\mod
2$.
 \pro{Theorem 3.1} Let $p$ be an odd prime, $a\in\Bbb Z_p$ and
$\ap\e 1\mod 2$. Then
$$\sum_{k=0}^{p-1}\b ak^2(-1)^k\e (-1)^{\f{p-\ap}2-1}
\f{2(a-\ap)}{a\b{p-\ap}{(p-\ap)/2}}\mod {p^2}.$$
\endpro
Proof. Using Lemma 3.1 we deduce that
$$\align &f_p(a)\e -\f{a+2}{4(a+1)}f_p(a+2)\e
\Big(-\f{a+2}{4(a+1)}\Big)\Big(-\f{a+4}{4(a+3)}\Big)f_p(a+4)
\e\cdots \\&\e \Big(-\f{a+2}{4(a+1)}\Big)\Big(-\f{a+4}{4(a+3)}\Big)
\cdots\Big(-\f{a+2(\f{p-\ap}2-1)}{4(a+2(\f{p-\ap}2-1)-1)}\Big)f_p(a+p-\ap-2)
\\&\e \f 1{(-16)^{\f{p-\ap}2-1}}\cdot\f{(a+1)(a+2)\cdots(a+p-\ap-2)}
{(\f{a+1}2)^2(\f{a+1}2+1)^2\cdots(\f{a+1}2+(\f{p-\ap}2-1)-1)^2}
\Big(-\f{a-\ap}2\Big)
\\&\e \f{(-1)^{\f{p-\ap}2}}{4^{p-\ap-2}}\cdot\f{(\ap+1)(\ap+2)\cdots(p-2)}{(-\f 32)^2(-\f 32-1)^2\cdots
(-\f 32-\f{p-\ap}2+2)^2}\cdot\f{a-\ap}2
\\&\e(-1)^{\f{p-\ap}2}4^{\ap}\f{(p-2)!/\ap!}{\b{-1/2}{(p-\ap)/2}^2\cdot
(\f{p-\ap}2)!^2}\cdot\f{a-\ap}2
\\&=(-1)^{\f{p-\ap}2}4^{\ap}\cdot\f{4^{p-\ap}\b{p-1}{\ap-1}\b{p-\ap}{(p-\ap)/2}}
{(p-1)\cdot\ap\cdot\b{p-\ap}{(p-\ap)/2}^2}\cdot\f {a-\ap}2
\\&\e -\f 2a(-1)^{\f{p-\ap}2}\f{a-\ap}{\b{p-\ap}{(p-\ap)/2}}\mod
{p^2}.
\endalign$$
This proves the theorem.

\pro{Theorem 3.2 } Let $p$ be an odd prime. Then
$$\align &\sum_{k=0}^{p-1}\f{\b{2k}k^2}{(-16)^k(2k-1)^2}
\\&\e \cases (-1)^{\f{p-1}4}\f px\mod{p^2}\q\t{if $p=x^2+4y^2\e 1\mod
4$ and $4\mid x-1$,}
\\(-1)^{\f{p+1}4}(2p+3-2^{p-1})\b{\f{p-1}2}{\f{p-3}4}
\mod {p^2}\q\t{if $p\e 3\mod 4$}
\endcases\endalign$$
and for $p\e 3\mod 4$,
$$\sum_{k=0}^{p-1}\f{\b{2k}k^2}{(-16)^k(2k-1)}
\e\f
12(-1)^{\f{p-3}4}(2p+3-2^{p-1})\b{\f{p-1}2}{\f{p-3}4}-(-1)^{\f{p-3}4}
\f p{\b{(p-1)/2}{(p-3)/4}}\mod {p^2}.$$

\endpro
Proof. Recall that $\b {1/2}k=-\b{2k}k\f 1{(-4)^k(2k-1)}$. For $p\e
1\mod 4$ and so $p=x^2+4y^2$ with $x\e 1\mod 4$, it is well known
that $\b{(p-1)/2}{(p-1)/4}\e 2x\mod p$ (see [BEW]). Now, taking
$a=\f 12$ and $\ap=\f{p+1}2$ in Theorem 3.1 gives
$$\sum_{k=0}^{p-1}\f{\b{2k}k^2}{(-16)^k(2k-1)^2}\e (-1)^{\f{p-1}4-1}
\f{2(\f 12-\f{p+1}2)}{\f 12\b{(p-1)/2}{(p-1)/4}}=(-1)^{\f{p-1}4}\f
px\mod {p^2}.$$ For $p\e 3\mod 4$, taking $a=\f 12$ and
$\ap=\f{p+1}2$ in (3.1) and then applying Remark 2.1 yields
$$\align \sum_{k=0}^{p-1}\f{\b{2k}k^2}{(-16)^k(2k-1)^2}&\e
(-1)^{\f{p+1}4}\b{\f{p+1}2}{\f{p+1}4}\big(1-\f
p2(H_{\f{p+1}2}-H_{\f{p+1}4}\big)\big) \\&\e(-1)^{\f{p+1}4}\cdot 2
\b{\f{p-1}2}{\f{p-3}4}\big(1-\f p2\big(\f 2{p+1}-2\qp 2-\big(\f
4{p+1}-3\qp 2\big)\big)\big)\\&=
(-1)^{\f{p+1}4}\b{\f{p-1}2}{\f{p-3}4}(2p+3-2^{p-1}) \mod
{p^2}.\endalign$$ By Theorem 2.1,
$$\align
&\sum_{k=0}^{p-1}\f{\b{2k}k^2}{(-16)^k(2k-1)^2}(2+4(2k-1)+4(2k-1)^2)
=\sum_{k=0}^{p-1}\b {1/2}k^2(-1)^k(1+(1-4k)^2) \\&\e
2p^2(1+2pH_{\f{p+1}2})\e 2p^2(1+2p(\f 2{p+1}-2\qp 2)\e
2p^2(1+4p(1-\qp 2)) \mod {p^4}.\endalign$$ By [Su3],
$$\sum_{k=0}^{p-1}\f{\b{2k}k^2}{(-16)^k}\e -\f{2(-1)^{\f{p+1}4}p}
{\b{(p+1)/2}{(p+1)/4}}=(-1)^{\f{p-3}4}\f p{\b{(p-1)/2}{(p-3)/4}}\mod
{p^2}.$$ Now, from the above we deduce that
$$\align&\sum_{k=0}^{p-1}\f{\b{2k}k^2}{(-16)^k(2k-1)}
\\&\e -\sum_{k=0}^{p-1}\f{\b{2k}k^2}{(-16)^k}-\f
12\sum_{k=0}^{p-1}\f{\b{2k}k^2}{(-16)^k(2k-1)^2}
\\&\e -(-1)^{\f{p-3}4}\f p{\b{(p-1)/2}{(p-3)/4}}-\f 12
(-1)^{\f{p+1}4}(2p+3-2^{p-1})\b{(p-1)/2}{(p-3)/4} \pmod {p^2}.
\endalign$$
This completes the proof.
\par{\bf Remark 3.1} Suppose that $p$ is an odd prime. The congruence for
$\sum_{k=0}^{p-1}\f{\b{2k}k^2}{(-16)^k(2k-1)^2}$ modulo $p^2$ was
conjectured by the author in [S8, Conjecture 5.48]. In [Su4], Z.W.
Sun obtained the congruence on
$\sum_{k=0}^{p-1}\f{\b{2k}k^2}{(-16)^k(2k-1)}$ modulo $p$ for $p\e
3\mod 4$, and modulo $p^2$ for $p\e 1\mod 4$.
\par\q
\par For any positive integers $n$ let
$$g_n(a)=\sum_{k=0}^{n-1}\b ak^2(-1)^k(a-2k).$$
For $k=0,1,2,\ldots$ set $F(a,k)=\b ak^2(-1)^k(a-2k)$ and
$$\align G(a,k)&=\f{a+2}{(a+1)^2}(-1)^{k-1}
\b{a+1}{k-1}^2
\\&\q\times\big(-4k^3+(14a+22)k^2-(a+2)(16a+20)k+(a+1)(5a^2+22a+23)\big).
\endalign$$
Using Maple one can easily check that
$$(a+1)F(a+2,k)+4(a+2)F(a,k)=G(a,k+1)-G(a,k).$$
Thus,
$$\align &(a+1)g_n(a+2)+4(a+2)g_n(a)
\\&=\sum_{k=0}^{n-1}\big((a+1)F(a+2,k)+4(a+2)F(a,k)\big)
=\sum_{k=0}^{n-1}(G(a,k+1)-G(a,k))\\&=G(a,n)-G(a,0)=G(a,n).
\endalign$$

\pro{Lemma 3.2} Suppose that $p$ is an odd prime and $a\in\Bbb Z_p$
with $a+1\not\e 0\mod p$. Then
$$g_p(a)\e \cases -\f{a+1}{4(a+2)}g_p(a+2)\mod {p^2}&\t{if $a+2\not\e 0\mod p$,}
\\\f 12+(a+2)\big(-\f 12+\qp 2\big)\mod{p^2}&\t{if $a+2\e 0\mod p$.}
\endcases$$
\endpro
Proof. For $\ap<p-2$ we see that $\b{a+1}{p-1}\e
\b{\ap+1}{p-1}=0\mod p$ and so $G(a,p)\e 0\mod {p^2}$. Hence,
$(a+1)g_p(a+2)+4(a+2)g_p(a)=G(a,p)\e 0\mod {p^2}$ and so the result
is true. For $\ap=p-2$ we  have $\b{a+1}{p-1}\e 1\mod {p^2}$ by the
proof of Lemma 3.1. Set $t=(a+2)/p$. We then have
$$\align \f{G(a,p)}{a+2}&=\f{1}{(a+1)^2}(-1)^{p-1}
\b{a+1}{p-1}^2
\\&\q\times\big(-4p^3+(14a+22)p^2-(a+2)(16a+20)p+(a+1)(5a^2+22a+23)\big)
\\&\e \f 1{pt-1}(5(pt-2)^2+22(pt-2)+23)\e\f{2pt-1}{pt-1}\e
-(a+1)\mod {p^2}.\endalign$$ Since $\b{p-1}r\e (-1)^r\mod p$ and
$$\sum_{k=1}^{p-1}\f{(-1)^k}k\e\sum_{k=1}^{p-1}\f{1+(-1)^k}k=
H_{\f{p-1}2}\e -2\qp 2\mod p,$$ we see that
$$\align g_p(a+2)&=pt+p^2t^2\sum_{k=1}^{p-1}\f 1{k^2}\b{pt-1}{k-1}^2(-1)^k(pt-2k)
\\&\e pt-2p^2t^2\sum_{k=1}^{p-1}\f{(-1)^k}k\e pt+4p^2t^2\qp 2\mod
{p^3}.\endalign$$ That is,
$$\f{g_p(a+2)}{a+2}\e 1+4(a+2)\qp 2\mod {p^2}.$$
Thus,
$$\align g_p(a)&=\f{G(a,p)}{4(a+2)}-\f{(a+1)g_p(a+2)}{4(a+2)}\e
-\f{a+1}4-\f{a+1}4(1+4(a+2)\qp 2)\\&\e -\f{a+1}2+(a+2)\qp 2=\f
12+(a+2)(-\f 12+\qp 2)\mod {p^2}.\endalign$$
 This completes the proof.

\pro{Theorem 3.3} Let $p$ be an odd prime, $a\in\Bbb Z_p$ and
$a\not\e 0,-1\mod p$. If $2\mid \ap$, then
$$\sum_{k=0}^{p-1}\b ak^2(-1)^k\Big(1-\f 2ak\Big)\e
(-1)^{\f{\ap}2}2^{2\ap}\f{a-\ap}{a\b{\ap}{\ap/2}} \mod {p^2}.$$ If
$2\nmid \ap$, then
$$\align\sum_{k=0}^{p-1}\b ak^2(-1)^k\Big(1-\f 2ak\Big)&\e
\f{a+1}a(-1)^{\f{p-\ap}2}\cdot 2^{2\ap+1}\b{p-2-\ap}{\f{p-2-\ap}2}
 \Big(1-2p\qp 2\\&\q+(p+a-\ap)(H_{\ap+1}-H_{\f{\ap+1}2})\Big)
\mod {p^2}.
\endalign$$
\endpro
 Proof.
 For $\ap\e 0\mod 2$, appealing to Lemma 3.2 we get
$$\align g_p(a)&\e -\f{4a}{a-1}g_p(a-2)\e
\Big(-\f{4a}{a-1}\Big)\Big(-\f{4(a-2)}{a-3}\Big) g_p(a-4)\e\cdots
\\&\e \Big(-\f{4a}{a-1}\Big)\Big(-\f{4(a-2)}{a-3}\Big)
\cdots \Big(-\f{4(a-\ap+2)}{a-\ap+1}\Big)g_p(a-\ap)
\\&=(-1)^{\f{\ap}2}\f{a(a-1)\cdots(a-\ap+1)}{(-\f{a-1}2)^2(-\f{a-3}2)^2\cdots
(-\f{a-\ap+1}2)^2}\\&\q\times\Big(a-\ap+\sum_{k=1}^{p-1}\f{(a-\ap)^2}{k^2}
\b{a-\ap-1}{k-1}^2(-1)^k(a-\ap-2k)\Big)
\\&\e (-1)^{\f{\ap}2}\f{\b a{\ap}\b
{\ap}{\ap/2}}{\b{-1/2
-(a-\ap)/2}{\ap/2}^2}(a-\ap)\e(-1)^{\f{\ap}2}\f{\b{\ap}{\ap/2}}{\b{-1/2}{\ap/2}^2}
(a-\ap)\\&=(-1)^{\f{\ap}2}2^{2\ap}\f{a-\ap}{\b{\ap}{\ap/2}} \mod
{p^2}.
\endalign$$
For $\ap\e 1\mod 2$, using Lemma 3.2 we see that
 $$\align &g_p(a)\e -\f{a+1}{4(a+2)}g_p(a+2)\e
 \Big(-\f{a+1}{4(a+2)}\Big)\Big(-\f{a+3}{4(a+4)}\Big)g_p(a+4)
 \e\cdots\\&\e \Big(-\f{a+1}{4(a+2)}\Big)\Big(-\f{a+3}{4(a+4)}\Big)
 \cdots\Big(-\f{a+2(\f{p-\ap}2-1)-1}{4(a+2(\f{p-\ap}2-1))}\Big)g_p\Big(a
 +2\big(\f{p-\ap}2-1\big)\Big)
 \\&=(-1)^{\f{p-\ap}2-1}\f{a+1}{p+a-\ap-1}\cdot
 \f{(\f{a+3}2)^2(\f{a+5}2)^2\cdots
 (\f{a+p-\ap-1}2)^2}
 {(a+2)(a+3)\cdots(a+(p-\ap-1))}\\&\q\times g_p(p+a-\ap-2)
\\&\e
\f{(-1)^{\f{p-\ap}2-1}(a+1)\b{\f{p+a-\ap}2-\f 12}{\f{p-\ap}2-1}^2 }
{(p+a-\ap-1)\b{p+a-\ap-1}{p-2-\ap}\b{p-2-\ap}{\f{p-2-\ap}2}}\\&\q\times
\Big(\f 12+(p+a-\ap)\big(-\f 12+\qp 2\big)\Big)\mod {p^2}.
\endalign$$
By [S7, (2.6)],
$$\align \b{p+a-\ap-1}{p-2-\ap}&\e
\b{-1}{p-2-\ap}\big(1+(p+a-\ap)(H_{p-1}-H_{\ap+1})\big)\\&\e
1-(p+a-\ap)H_{\ap+1}\mod {p^2}\endalign$$ and
$$\align&\b{\f{p+a-\ap}2-\f 12}{\f{p-\ap}2-1}\e\b{-\f
12}{\f{p-\ap}2-1}\Big(1+\f{p+a-\ap}2(H_{\f{p-1}2}-H_{\f{\ap+1}2})\Big)
\\&\e (-4)^{-(\f{p-\ap}2-1)}\b{p-\ap-2}{\f{p-\ap}2-1}
\Big(1-\f{p+a-\ap}2(2\qp 2+H_{\f{\ap+1}2})\Big)\mod{p^2}.\endalign$$
Hence,
$$\align &\f{\b{\f{p+a-\ap}2-\f 12}{\f{p-\ap}2-1}^2}
{\b{p+a-\ap-1}{p-2-\ap}\b{p-2-\ap}{\f{p-2-\ap}2}} \e
4^{-(p-2-\ap)}\b{p-2-\ap}{\f{p-2-\ap}2} \f{(1-\f{p+a-\ap}2(2\qp
2+H_{\f{\ap+1}2}))^2}{1-(p+a-\ap)H_{\ap+1}}
\\&\e 4^{\ap+1}(1-2p\qp 2)\b{p-2-\ap}{\f{p-2-\ap}2}
\\&\q\times(1-(p+a-\ap)(2\qp 2+H_{\f{\ap+1}2}))(1+(p+a-\ap)H_{\ap+1})
\\&\e  4^{\ap+1}\b{p-2-\ap}{\f{p-2-\ap}2}
\\&\q\times\big(1-2p\qp 2+(p+a-\ap)(-2\qp 2+H_{\ap+1}-H_{\f{\ap+1}2})\big) \mod
{p^2}.\endalign$$ It then follows that
$$\align g_p(a)&\e
(-1)^{\f{p-\ap}2}\f{a+1}{1-(p+a-\ap)}\cdot4^{\ap+1}\b{p-2-\ap}{\f{p-2-\ap}2}
\\&\q\times\big(1-2p\qp 2+(p+a-\ap)(-2\qp 2+H_{\ap+1}-H_{\f{\ap+1}2})\big)
\\&\q\times \f 12\big(1+(p+a-\ap)(-1+2\qp 2)\big)
\\&\e (-1)^{\f{p-\ap}2}(a+1)\cdot 2^{2\ap+1}\b{p-2-\ap}{\f{p-2-\ap}2}
\\&\q\times \big(1-2p\qp 2+(p+a-\ap)(H_{\ap+1}-H_{\f{\ap+1}2})\big)
\mod {p^2}.
\endalign$$
This completes the proof.

\pro{Corollary 3.1} Let $p>3$ be a prime. Then
$$\sum_{k=0}^{p-1}\b{-\f 13}k^2(-1)^k(6k+1)
\e\cases \f{(-1)^{\f{p-1}2}p}{2^{\f{p-1}3}\b{(p-1)/3}{(p-1)/6}}\mod
{p^2}&\t{if $3\mid p-1$,}
\\(-1)^{\f{p-5}6}2^{\f{p+4}3}\b{\f{p-5}3}{\f{p-5}6}\big(4-
p-2^{p}\big)\mod {p^2}&\t{if $3\mid p-2$.}
\endcases$$
\endpro
Proof. Taking $a=-\f 13$ in Theorem 3.3 and then applying the fact
$H_{p-1-k}\e H_k\mod p$ yields the result.

\par{\bf Remark 3.2} For any positive integers $n$ we have the
identities:
$$\align &\sum_{k=0}^{n-1}\b ak^2\Big(1-\f
2ak\Big)^3=\f{4n^2-4(a+1)n+a(a+3)}{a(a-1)}\b{a-1}{n-1}^2,
\\&\sum_{k=0}^{n-1}\b ak^2\Big(1-\f
2ak\Big)^5=\f{n^2\b
an^2}{a^5(a-1)(a-2)}\big(16(a-1)n^4-32(a^2-1)n^3\\&\qq\qq\qq\qq\q+
(24a^3+32a^2-48a-16)n^2+(-8a^4-24a^3+16a^2+32a)n\\&\qq\qq\qq\qq\q
+a^5+5a^4+2a^3-16a^2\big).
\endalign$$
Thus, for any odd prime $p$ and $a\in\Bbb Z_p$ with $a\not\e
0,1,2\mod p$,
$$\sum_{k=0}^{p-1}\b ak^2\Big(1-\f
2ak\Big)^3\e \sum_{k=0}^{p-1}\b ak^2\Big(1-\f 2ak\Big)^5\e 0\mod
{p^2}.$$ Furthermore, we conjecture that for any odd prime $p$,
$r\in\{0,1,2,\ldots\}$ and $a\in\Bbb Z_p$ with $\ap>r$,
$$\sum_{k=0}^{p-1}\b ak^2\Big(1-\f
2ak\Big)^{2r+1}\e 0\mod {p^2}.$$

\section*{4. Congruences for $\sum_{k=0}^{p-1}(1-\f 2ak)\b ak^3$
and $\sum_{k=0}^{p-1}(1-\f 2ak)^3\b ak^3$
 modulo $p^4$}
\par
 For $n=1,2,3,\ldots$ let
 $$C_n(a)=\sum_{k=0}^{n-1}(a-2k)\b ak^3.$$
Set $G(a,k)=\f{(2a+2-k)k^3}{(a+1-k)^3}\b ak^3$. It is easy to check
that
$$(a-2k)\b ak^3+(a+1-2k)\b{a+1}k^3=G(a,k+1)-G(a,k).$$
Thus,
$$\sum_{k=0}^{n-1}\Big((a-2k)\b ak^3+(a+1-2k)\b{a+1}k^3\Big)
=\sum_{k=0}^{n-1}(G(a,k+1)-G(a,k))=G(a,n).$$ Hence,
$$C_n(a)+C_n(a+1)=\f{(2a+2-n)n^3}{(a+1-n)^3}\b an^3.\tag 4.1$$

\pro{Lemma 4.1} Let $p$ be an odd prime, $a\in\Bbb Z_p$ and $a\not\e
-1\mod p$. Then
$$C_p(a)+C_p(a+1)\e \f{2(a-\ap)^3}{(\ap+1)^2}\mod {p^4}.$$
\endpro
Proof. Set $t=(a-\ap)/p$. By Lucas' theorem, $\b
ap=\b{\ap+pt}{0+1\cdot p}\e t\mod p$. Now, taking $n=p$ in (4.1) we
see that
$$C_p(a)+C_p(a+1)=\f{(2a+2-p)p^3}{(a+1-p)^3}\b ap^3
\e \f 2{(a+1)^2}p^3\cdot t^3\e\f {2(a-\ap)^3}{(\ap+1)^2}\mod
{p^4}.$$ This proves the lemma.

\pro{Theorem 4.1} Let $p$ be an odd prime, $a\in\Bbb Z_p$ and
$a\not\e 0\mod p$. Then
$$\align&\sum_{k=0}^{p-1}\Big(1-\f 2ak\Big)\b ak^3\e
(-1)^{\ap}\f{a-\ap}{a}+\f{(a-\ap)^3}aE_{p-3}(-a)\mod {p^4},\tag 4.2
\\&\sum_{k=0}^{p-1}\Big(1-\f 2ak\Big)^3\b ak^3
\e -3(-1)^{\ap}\f{a-\ap}{a}\tag
4.3\\&\qq\qq\qq\qq\qq+\f{(a-\ap)^3}{a^3}\big(4-3a^2E_{p-3}(-a)\big)\mod
{p^4}.\endalign$$
\endpro
Proof. By [S6, Lemma 2.2],
$$\sum_{k=1}^{\ap}\f{(-1)^k}{k^2}\e \f 12(-1)^{\ap}E_{p-3}(-a)\mod
p.$$For $k=1,2,\ldots,\ap$ we see that $a-k-\langle
a-k\rangle_p=a-k-(\ap-k)=a-\ap$ and so $C_p(a-k)+C_p(a-k+1)\e
\f{2(a-\ap)^3}{(\ap-k+1)^2}\mod {p^4}$. Thus,
$$\align
&C_p(a)+(-1)^{\ap-1}C_p(a-\ap)\\&=\sum_{k=1}^{\ap}(-1)^{k-1}\big(C_p(a-k)+C_p(a-k+1)\big)
\\&\e \sum_{k=1}^{\ap}(-1)^{k-1}\f{2(a-\ap)^3}{(\ap-k+1)^2}
=2(-1)^{\ap}(a-\ap)^3\sum_{r=1}^{\ap}\f{(-1)^r}{r^2}\\&\e
(a-\ap)^3E_{p-3}(-a)\mod {p^4}.\endalign$$  Since
$$\sum_{k=1}^{p-1}\f{(-1)^k}{k^2}=\f 12\sum_{r=1}^{(p-1)/2}\f
1{r^2}-\sum_{k=1}^{p-1}\f 1{k^2}\e 0-0=0\mod p,$$ we have
$$\align
C_p(a-\ap)&=a-\ap+\sum_{k=1}^{p-1}(a-\ap-2k)\f{(a-\ap)^3}{k^3}\b{a-\ap-1}{k-1}^3
\\&\e a-\ap+\sum_{k=1}^{p-1}(a-\ap-2k)\f{(a-\ap)^3}{k^3}(-1)^{k-1}
\\&\e a-\ap+2(a-\ap)^3\sum_{k=1}^{p-1}\f{(-1)^k}{k^2}\e a-\ap\mod {p^4}.\endalign$$
Now, combining all the above yields (4.2). By Theorem 2.1,
$$\sum_{k=0}^{p-1}\b ak^3\Big(3\Big(1-\f 2ak\Big)+\Big(1-\f
2ak\Big)^3\Big)\e 4\f{(a-\ap)^3}{a^3}\mod {p^4}.$$ This together
with (4.2) yields (4.3).

\pro{Theorem 4.2} Let $p>3$ be a prime. Then
$$\align
&\sum_{k=0}^{p-1}(4k+1)^3\f{\b{2k}k^3}{(-64)^k}\e
-3(-1)^{\f{p-1}2}p+p^3(4-3E_{p-3})\mod {p^4},
\\&\sum_{k=0}^{p-1}(4k-1)\f{\b{2k}k^3}{(-64)^k(2k-1)^3}\e
(-1)^{\f{p-1}2}p+p^3(E_{p-3}-2)\mod {p^4},
\\&\sum_{k=0}^{p-1}(4k-1)^3\f{\b{2k}k^3}{(-64)^k(2k-1)^3}\e
-3(-1)^{\f{p-1}2}p+p^3(2-3E_{p-3})\mod {p^4},
\\&\sum_{k=0}^{p-1}(12k+1)\b{-\f 16}k^3
\e\cases (-1)^{\f{p-1}2}p+\f 59p^3E_{p-3}\mod {p^4}&\t{if $3\mid
p-1$,}
\\5(-1)^{\f{p-1}2}p+\f{625}9p^3E_{p-3}\mod {p^4}&\t{if $3\mid p-2$,}
\endcases
\\&\sum_{k=0}^{p-1}(12k+1)^3\b{-\f 16}k^3
\e\cases -3(-1)^{\f{p-1}2}p+p^3\big(4-\f 53E_{p-3}\big)\mod
{p^4}&\t{if $3\mid p-1$,}
\\-15(-1)^{\f{p-1}2}p+p^3\big(500-\f{625}3E_{p-3}\big)\mod {p^4}
&\t{if $3\mid p-2$,}
\endcases
\\&\sum_{k=0}^{p-1}(6k+1)\b{-\f 13}k^3
\e\cases p+p^3U_{p-3}\mod {p^4}&\t{if $3\mid p-1$,}
\\-2p+8p^3U_{p-3}\mod {p^4}&\t{if $3\mid p-2$,}\endcases
\\&\sum_{k=0}^{p-1}(6k+1)^3\b{-\f 13}k^3
\e\cases -3p+p^3(4-3U_{p-3})\mod {p^4}&\t{if $3\mid p-1$,}
\\6p+8p^3(4-3U_{p-3})\mod {p^4}&\t{if $3\mid p-2$.}\endcases
\endalign$$
\endpro
Proof. From [S4, pp.3300-3301] and [MOS] we know that
$$\align &E_{p-3}\Ls 12\e 4E_{p-3}\mod p,\q E_{p-3}\Big(-\f 12\Big)=\f
1{2^{p-4}}-E_{p-3}\Ls 12\e 8-4E_{p-3}\mod p,
\\&E_{p-3}\Ls 13\e 9U_{p-3}\mod p,\q
E_{p-3}\Ls 16\e 20E_{p-3}\mod p.\endalign$$ Recall that
$\b{-1/2}k=\b{2k}k\f 1{(-4)^k}$ and $\b {1/2}k=-\b{2k}k\f
1{(-4)^k(2k-1)}$. Taking $a=-\f 12,\f 12,-\f 13,-\f 16$ in Theorem
4.1 and then applying the above yields the result.
\par\q
\par{\bf Remark 4.1} Let $p$ be an odd prime, $r\in\{0,1,2,\ldots\}$
 and $a\in\Bbb Z_p$ with
$a\not\e 0\mod p$. Using the WZ method, one may similarly deduce the
congruences for $\sum_{k=0}^{p-1}\b ak^3(1-\f 2ak)^{2r+1}$ and
$\sum_{k=0}^{p-1}\b ak^3(-1)^k(1-\f 2ak)^r$
 modulo $p^3$.

\section*{5. The congruence for $\sum_{k=0}^{p-1}(a-2k)\b ak^4$ modulo $p^4$}
\par For $n=1,2,3,\ldots$ let
$$Q_n(a)=\sum_{k=0}^{n-1}(a-2k)\b{a}k^4.$$
 For $k=0,1,2,\ldots$ set
$F(a,k)=(a-2k)\b ak^4$ and
$$G(a,k)=\f{k^4(2k^2-6(a+1)k+5(a+1)^2)}{(a+1-k)^4}\b ak^4.$$
It is easy to check that
$$(a+1)F(a+1,k)+2(2a+1)F(a,k)=G(a,k+1)-G(a,k).$$
Thus,
$$\align &(a+1)Q_n(a+1)+2(2a+1)Q_n(a)
\\&=(a+1)\sum_{k=0}^{n-1}F(a+1,k)+2(2a+1)\sum_{k=0}^{n-1}F(a,k)
\\&=\sum_{k=0}^{n-1}(G(a,k+1)-G(a,k))=G(a,n)-G(a,0)=G(a,n).\endalign$$
That is,
$$\aligned &(a+1)Q_n(a+1)+2(2a+1)Q_n(a)
\\&=\f{2n^2-6(a+1)n+5(a+1)^2}{(a+1-n)^4}\cdot n^4\b an^4.\endaligned\tag 5.1$$
\pro{Lemma 5.1} Let $p$ be an odd prime, $a\in\Bbb Z_p$ and $a\not\e
-1\mod p$. Then
$$(a+1)Q_p(a+1)\e -2(2a+1)Q_p(a)\mod {p^4}.$$
\endpro
Proof. Taking $n=p$ in (5.1) yields the result.

\pro{Theorem 5.1} Let $p$ be an odd prime, $a\in\Bbb Z_p$, $a\not\e
0,-1\mod p$ and $a'=(a-\ap)/p$. If $\ap<\f p2$, then
$$\align Q_p(a)&\e a'p(-1)^{\ap}\b{2\ap}{\ap}\Big(1+2a'p(H_{2\ap}-H_{\ap})
\\&\qq+{a'}^2p^2\big(2(H_{2\ap}-H_{\ap})^2-2H_{2\ap}^{(2)}+H_{\ap}^{(2)}\big)
\Big)\mod {p^4}.\endalign$$ If $\ap>\f p2$, then
$$\align Q_p(a)&\e \f{a'(2a'+1)}{2a+1}p^2(-1)^{\f{p-1}2}(-16)^{\ap}
\b{2(\ap-\f{p-1}2)}{\ap-\f{p-1}2}^{-1}
\Big(1+p\Big(\f{2a'+1}{2a+1}\\&\q+(2a'-1)\qp
2-\f{2a'+1}2\big(H_{p-1-\ap}-H_{\ap-\f{p+1}2}\big)\Big)\Big)
 \mod {p^4}.\endalign$$
\endpro
Proof. By Lemma 5.1,
$$\align Q_p(a)&\e \f{-2(2a-1)}aQ_p(a-1)\e
\f{-2(2a-1)}a\cdot\f{-2(2(a-1)-1)}{a-1}Q_p(a-2)
\\&\e\cdots\e\prod_{k=0}^{\ap-1}\f{-2(2(a-k)-1)}{a-k}\cdot
Q_p(a-\ap)\\&=(-4)^{\ap}\prod_{k=0}^{\ap-1}\f{a-\f 12-k}{a-k}\cdot
(a-\ap)\mod {p^4},\endalign$$ which implies
 $$Q_p(a)\e(-4)^{\ap}(a-\ap)\f{\b{a-1/2}{\ap}}{\b a{\ap}}\mod {p^4}.
\tag 5.2$$
 For $1\le n\le m\le p-1$ and $t\in\Bbb Z_p$ we
see that
$$\align &\b{m+pt}n=\f{(pt+m)(pt+m-1)\cdots(pt+m-(n-1))}{n!}\\&\e
\b mn\Big(1+pt\sum_{m-n+1\le k\le m}\f 1k+p^2t^2\sum_{m-n+1\le
i<j\le m}\f 1{ij}\Big)
\\&=\b mn\Big(1+pt(H_m-H_{m-n})+\f
12p^2t^2\big((H_m-H_{m-n})^2-(H_m^{(2)}-H_{m-n}^{(2)})\big)\Big)
\mod {p^3}.
\endalign$$
Thus,
$$\b a{\ap}=\b{\ap+a'p}{\ap}
\e 1+a'pH_{\ap}+\f 12{a'}^2p^2(H_{\ap}^2-H_{\ap}^{(2)})\mod {p^3}.$$
For $1\le \ap\le \f{p-1}2$ we see that
$$\align &\b{a-1/2}{\ap}\\&
=(-1)^{\ap}\b{-a'p-\f 12}{\ap} =(-1)^{\ap}\f{(-a'p-\f 12)(-a'p-\f
32)\cdots(-a'p-\f{2\ap-1}2)}{\ap!}
\\&\e (-1)^{\ap}\f{-\f 12(-\f 12-1)\cdots(-\f 12-\ap+1)}{\ap !}
\\&\q\times\Big(1-a'p\sum_{r=1}^{\ap}\f 1{-\f{2r-1}2}+{a'}^2p^2\sum_{1\le
i<j\le \ap}\f 1{(-\f{2i-1}2)(-\f{2j-1}2)}\Big)
\\&=(-1)^{\ap}\b{-1/2}{\ap}\Big(1+2a'p\sum_{r=1}^{\ap}\f 1{2r-1}
+4{a'}^2p^2\sum_{1\le i<j\le \ap}\f 1{(2i-1)(2j-1)}\Big)
\\&=\f {\b{2\ap}{\ap}}{4^{\ap}}\Big(1+2a'p\sum_{r=1}^{\ap}\f 1{2r-1}
+2{a'}^2p^2\Big(\big(\sum_{r=1}^{\ap}\f
1{2r-1}\big)^2-\sum_{r=1}^{\ap}\f 1{(2r-1)^2}\Big)\Big) \mod {p^3}.
\endalign$$
Since $(1+bp+cp^2)(1-bp+(b^2-c)p^2)\e 1\mod {p^3}$ for $b,c\in\Bbb
Z_p$ and
$$\sum_{r=1}^{\ap}\f
1{(2r-1)^i}=\sum_{k=1}^{2\ap}\f 1{k^i}-\sum_{k=1}^{\ap}\f 1{(2k)^i}
=H_{2\ap}^{(i)}-\f 1{2^i}H_{\ap}^{(i)}\qtq{for}i=1,2,$$ from the
above we deduce that for $1\le \ap\le\f{p-1}2$,
$$\align Q_p(a)&\e (-4)^{\ap}(a-\ap)\f{\b{a-1/2}{\ap}}{\b a{\ap}}
\\&\e (-1)^{\ap}a'p\b{2\ap}{\ap}\f{1+2a'p\sum_{r=1}^{\ap}\f 1{2r-1}
+2{a'}^2p^2\big(\big(\sum_{r=1}^{\ap}\f
1{2r-1}\big)^2-\sum_{r=1}^{\ap}\f 1{(2r-1)^2}\big)}{1+a'pH_{\ap}+\f
12{a'}^2p^2(H_{\ap}^2-H_{\ap}^{(2)})}
\\&\e (-1)^{\ap}a'p\b{2\ap}{\ap}\Big(1+a'p(2H_{2\ap}-H_{\ap})
\\&\q+2{a'}^2p^2\big(\big(H_{2\ap}-\f
12H_{\ap}\big)^2-\big(H_{2\ap}^{(2)}-\f 14H_{\ap}^{(2)}\big)\Big)
\\&\q\times\Big(1-a'pH_{\ap}+\f
12{a'}^2p^2(H_{\ap}^2+H_{\ap}^{(2)})\Big)
\\&\e (-1)^{\ap}a'p\b{2\ap}{\ap}\Big(1+2a'p(H_{2\ap}-H_{\ap})
\\&\qq+{a'}^2p^2\big(2(H_{2\ap}-H_{\ap})^2-2H_{2\ap}^{(2)}+H_{\ap}^{(2)}\big)
\Big)\mod {p^4}.\endalign$$
For $\ap\ge \f{p+1}2$ and $s=a'+\f 12$
 we see that
 $$\align \b{a-\f 12}{\ap}&=\b{\ap-\f{p+1}2+sp}{\ap}
 =\f{(\ap-\f{p+1}2+sp)\cdots
 (1+sp)sp(-1+sp)\cdots(-\f{p-1}2+sp)}{\ap!}
 \\&\e sp\cdot\f{\big(\ap-\f{p+1}2\big)!(1+spH_{\ap-\f{p+1}2})
 \cdot (-1)^{\f{p-1}2}\ls{p-1}2!\big(1-spH_{\f{p-1}2}\big)}{\ap!}
 \\&\e (-1)^{\f{p-1}2}sp\big(1-sp(H_{\f{p-1}2}-H_{\ap-\f{p+1}2})\big)
 \cdot\f 1{(\ap-\f{p-1}2)\b{\ap}{(p-1)/2}}\mod {p^3}.
 \endalign$$
 By Remark 2.1, $H_{p-1-k}\e H_k\mod p$ for
$k\in\{0,1,\ldots,p-1\}$ and $H_{\f{p-1}2}\e -2\qp 2\mod p$. By
[S7], $2H_{2k}-H_k\e H_{\f{p-1}2+k}-H_{\f{p-1}2}\mod p$ for $1\le
k\le \f{p-1}2$. Thus,
$$\align \b{\ap}{\f{p-1}2}&=(-1)^{\ap-\f{p-1}2}\b{-\f 12-\f
p2}{\ap-\f{p-1}2}\\&\e (-1)^{\ap-\f{p-1}2}\b{-\f 12}{\ap-\f{p-1}2}
\Big(1-\f p2\sum_{r=1}^{\ap-\f{p-1}2}\f 1{-(2r-1)/2}\Big)
\\&=\f{\b{2(\ap-\f{p-1}2)}{\ap-\f{p-1}2}}{4^{\ap-\f{p-1}2}}
\Big(1+p(H_{2(\ap-\f{p-1}2)}-\f 12H_{\ap-\f{p-1}2})\Big)
\\&\e 2^{p-1}\f{\b{2(\ap-\f{p-1}2)}{\ap-\f{p-1}2}}{4^{\ap}}
\Big(1+\f p2(H_{\ap}-H_{\f{p-1}2})\Big)
\\&\e\f{\b{2(\ap-\f{p-1}2)}{\ap-\f{p-1}2}}{4^{\ap}}(1+p\qp 2)
\Big(1+\f p2(H_{p-1-\ap}+2\qp 2)\Big)
\\&\e \f{\b{2(\ap-\f{p-1}2)}{\ap-\f{p-1}2}}{4^{\ap}}\Big(1+p\Big(2\qp 2+\f
12H_{p-1-\ap}\Big)\Big)\mod {p^2}
\endalign$$
 and so
$$\align Q_p(a)&\e (-4)^{\ap}a'p\cdot(-1)^{\f{p-1}2}\Big(a'+\f
12\Big)p \f{1-(a'+\f
12)p(H_{\f{p-1}2}-H_{\ap-\f{p+1}2})}{(\ap-\f{p-1}2)\b{\ap}{(p-1)/2}(1+a'pH_{\ap})}
\\&\e (-1)^{\f{p-1}2}(-4)^{\ap}a'\Big(a'+\f
12\Big)p^2\\&\q\times\Big(1-p\big(a'H_{\ap}+\big(a'+\f
12\big)\big(H_{\f{p-1}2}-H_{\ap-\f{p+1}2}\big)\big)\Big) \f 1{(a+\f
12-(a'+\f 12)p)\b{\ap}{(p-1)/2}}
\\&\e (-1)^{\f{p-1}2}p^2(-4)^{\ap}a'\Big(a'+\f
12\Big)\f{a+\f 12+(a'+\f 12)p}{(a+\f
12)^2}\\&\q\times\Big(1-p\Big(a'H_{p-1-\ap}-\big(a'+\f
12\big)\big(2\qp 2+H_{\ap-\f{p+1}2}\big)\Big)\Big)
\\&\q\times4^{\ap}\b{2(\ap-\f{p-1}2)}{\ap-\f{p-1}2}^{-1}\Big(1-p\Big(2\qp 2+\f
12H_{p-1-\ap}\Big)\Big)\Big)
\\&\e  (-1)^{\f{p-1}2}p^2(-16)^{\ap}
\f{a'(2a'+1)}{2a+1}\b{2(\ap-\f{p-1}2)}{\ap-\f{p-1}2}^{-1}
\Big(1+p\Big(\f{2a'+1}{2a+1}\\&\q+(2a'-1)\qp 2-\Big(a'+\f
12\Big)\big(H_{p-1-\ap}-H_{\ap-\f{p+1}2}\big)\Big)\Big)
 \mod {p^4}.
\endalign$$
 This completes the proof.

 \pro{Corollary 5.1 ([L, Theorem 1.1])} Let $p>3$ be a prime. Then
$$\sum_{k=0}^{(p-1)/2}(4k+1)\f{\b{2k}k^4}{256^k}\e p\mod {p^4}.$$
\endpro
Proof. Taking $a=-\f 12$ and $\ap=\f{p-1}2$ in (5.2) and then
applying Morley's congruence $(-1)^{\f{p-1}2}\b{p-1}{(p-1)/2}\e
4^{p-1}\mod {p^3}$ gives
$$Q_p\Big(-\f 12\Big)\e (-4)^{\f{p-1}2}\Big(-\f
p2\Big)\f{\b{-1}{(p-1)/2}}{\b{-1/2}{(p-1)/2}}=-\f
p2(-1)^{\f{p-1}2}\f{(-4)^{p-1}}{\b{p-1}{(p-1)/2}}\e -\f p2\mod
{p^4}.$$ Hence,
$$\sum_{k=0}^{(p-1)/2}(4k+1)\f{\b{2k}k^4}{256^k}\e
\sum_{k=0}^{p-1}(4k+1)\b{1/2}k^4=-2Q_p\Big(-\f 12\Big)\e p\mod
{p^4}.$$
\par By
calculations with Maple, we make the following conjecture.
 \pro{Conjecture 5.1} Let $p>3$ be a prime. Then
$$\sum_{k=0}^{p-1}(4k+1)\f{\b{2k}k^4}{256^k}\e
\sum_{k=0}^{(p-1)/2}(4k+1)\f{\b{2k}k^4}{256^k}\e p+\f
76p^4B_{p-3}\mod {p^5}.$$
\endpro
\pro{Corollary 5.2} Let $p$ be an odd prime. Then
$$\sum_{k=0}^{(p-1)/2}(4k+3)\f{C_k^4}{256^k}\e 16\mod {p^4},$$
where $C_n$ is the Catalan number given by $C_n=\b{2n}n\f 1{n+1}$.
\endpro
Proof. Set $a=\f 12$. Then $\ap=\f{p+1}2>\f p2$ and $a'=-\f 12$. By
Theorem 5.1,
$$\sum_{k=0}^{p-1}\Big(\f 12-2k\Big)\b{1/2}k^4\e 0\mod {p^4}.$$
Observe that $p\mid\b{2k}k$ for $\f p2<k<p$ and
$$\b{1/2}k=\f
1{1-2k}\b{-1/2}k=-\f{\b{2k}k}{(2k-1)(-4)^k}=-\f{2C_{k-1}}{(-4)^k}
\qtq{for}k\ge 1.$$ We then have
$$\align &\sum_{k=0}^{(p-1)/2}(4k+3)\f{C_k^4}{256^k}
\e \sum_{r=0}^{p-2}(4r+3)\f{C_r^4}{256^r}
=\sum_{k=1}^{p-1}(4k-1)\f{C_{k-1}^4}{256^{k-1}}
\\&=\sum_{k=1}^{p-1}(4k-1)\f{\b{1/2}k^4\cdot (-4)^{4k}}{16\cdot
256^{k-1}}=16+16\sum_{k=0}^{p-1}(4k-1)\b{1/2}k^4\e 16\mod{p^4}.
\endalign$$
This proves the corollary.
\par{\bf Remark 5.1} Let $p$ be an odd prime. From [LW, Corollary 1.2]
and Theorem 2.1 (with $a=\f 12$) one may prove the following
stronger congruence:
$$\sum_{k=0}^{p-1}(4k+3)\f{C_k^4}{256^k}\e
\sum_{k=0}^{(p-1)/2}(4k+3)\f{C_k^4}{256^k}\e 16+80p^4\mod {p^5}.$$

\par For any positive integers $n$ let
$$S_n(x)=\sum_{k=0}^{n-1}\b{2k}k\b xk\b{-1-x}k\f 1{4^k}.$$
 \pro{Theorem 5.2} Let $p$ be a prime, $a\in\Bbb Z$, $1\le
\ap<\f p2$ and $a'=(a-\ap)/p$. Then
$$\Ls{Q_p(a)}p^2\e 2^{4\ap}{a'}^2S_p(2a)
\big(1+4a'p\qp 2+2a'(4a'-1)p^2\qp 2^2\big)\mod {p^3}.$$
\endpro
Proof. By [S9, Theorem 2.1],
$$S_p(2a)\e \f{\b{a-1/2}{\ap}^2}{\b a{\ap}^2}\big(1-4a'p\qp
2+2a'(4a'+1)p^2\qp 2^2\big)\mod {p^3}.$$ By (5.2),
$$\f{Q_p(a)}p\e(-4)^{\ap}a'\f{\b{a-1/2}{\ap}}{\b a{\ap}}\mod
{p^3}.$$ Recall that $(1+bp+cp^2)(1-bp+(b^2-c)p^2)\e 1\mod {p^3}$
for $b,c\in\Bbb Z_p$. We then have
$$\align \Ls{Q_p(a)}p^2&\e 16^{\ap}{a'}^2\f{\b{a-1/2}{\ap}^2}{\b a{\ap}^2}
\e 16^{\ap}{a'}^2\f{S_p(2a)}{1-4a'p\qp 2+2a'(4a'+1)p^2\qp 2^2}
\\&\e 2^{4\ap}{a'}^2S_p(2a)\big(1+4a'p\qp 2+2a'(4a'-1)p^2\qp 2^2\big)\mod {p^3}.
\endalign$$
This completes the proof.
\par\q
\pro{Lemma 5.2 ([S1, Theorem 5.2] and [S2, Corollaries 3.3 and
3.7])} Let $p>3$ be a prime. Then
$$\align &H_{\f{p-1}2}\e -2\qp 2+p\qp 2^2\mod {p^2},
\q H_{\f{p-1}2}^{(2)}\e 0\mod p,
\\&H_{[\f p4]}\e -3\qp 2+\f 32p\qp 2^2-(-1)^{\f{p-1}2}pE_{p-3}\mod
{p^2},
\\&H_{[\f p4]}^{(2)}\e 4(-1)^{\f{p-1}2}E_{p-3}\mod p.
\endalign$$

\pro{Theorem 5.3} Let $p>3$ be a prime. Then
$$\align &\sum_{k=0}^{p-1}(8k+1)\b{-1/4}k^4
\\&\e\cases (-1)^{\f{p-1}4}\big(2xp-\f{p^2}{2x}-\f{p^3}{8x^3}\big)
\mod {p^4}&\t{if $p=x^2+y^2\e 1\mod 4$ and $4\mid x-1$,}
\\3(-1)^{\f{p+1}4}\f{(2p-1-2^{p-1})p^2}{\binom{(p-1)/2}{(p-3)/4}}
\mod {p^4}&\t{if $p\e 3\mod 4$.}
\endcases\endalign$$
\endpro
Proof. We first assume that $p\e 1\mod 4$ and so $p=x^2+y^2$ with
$4\mid x-1$. Set $a=-\f 14$. Then $\ap=\f{p-1}4$ and $a'=-\f 14$. By
Theorem 5.1,
$$\align Q_p\Big(-\f 14\Big)&
\e (-1)^{\f{p-1}4}\Big(-\f p4\Big)\b{(p-1)/2}{(p-1)/4} \Big(1-\f
p2(H_{\f{p-1}2}-H_{\f{p-1}4})
\\&\q+\f{p^2}{16}\big(2(H_{\f{p-1}2}-H_{\f{p-1}4})^2-2H_{\f{p-1}2}^{(2)}
+H_{\f{p-1}4}^{(2)}\big)\Big) \mod {p^4}.
\endalign$$
From Lemma 5.2 we have
$$H_{\f{p-1}2}-H_{[\f p4]}\e \qp 2-\f p2\qp 2^2+(-1)^{\f{p-1}2}pE_{p-3}\mod {p^3}.
\tag 5.3$$ By [CD],
$$\b{\f{p-1}2}{\f{p-1}4}\e
\Big(2x-\f p{2x}-\f{p^2}{8x^3}\Big)\Big(1+\f p2\qp
2+\f{p^2}8(2E_{p-3}-\qp 2^2)\Big)\mod {p^3}.\tag 5.4$$ Now, from the
above and Lemma 5.2 we deduce that
$$\align \sum_{k=0}^{p-1}(8k+1)\b{-1/4}k^4&=-4Q_p\Big(-\f 14\Big)\e
(-1)^{\f{p-1}4}p\b{(p-1)/2}{(p-1)/4}\\&\q\times\Big(1-\f p2\Big(\qp
2-\f p2\qp 2^2+pE_{p-3}\Big)+\f{p^2}{16}\big(2\qp
2^2+4E_{p-3}\big)\Big)
\\&\e
(-1)^{\f{p-1}4}p\b{(p-1)/2}{(p-1)/4}\Big(1-\f 12p\qp 2+p^2\Big(\f
38\qp 2^2-\f 14E_{p-3}\Big)\Big)
\\&\e (-1)^{\f{p-1}4}p\Big(2x-\f p{2x}-\f{p^2}{8x^3}\Big)\Big(1+\f p2\qp
2+\f{p^2}8(2E_{p-3}-\qp 2^2)\Big)\\&\q\times\Big(1-\f 12p\qp
2+p^2\Big(\f 38\qp 2^2-\f 14E_{p-3}\Big)\Big)
\\&\e (-1)^{\f{p-1}4}p\Big(2x-\f p{2x}-\f{p^2}{8x^3}\Big)\mod {p^4}.
\endalign$$
\par Suppose $p\e 3\mod 4$ and $a=-\f 14$. Then $\ap=\f{3p-1}4$
and $a'=-\f 34$. By Theorem 5.1,
$$\align \sum_{k=0}^{p-1}(8k+1)\b{-\f 14}k^4&=-4Q_p\Big(-\f 14\Big)\e
4p^2(-1)^{\f{3p-1}4}2^{3p-1}\f{(-\f 34)(-\f
12)}{\b{(p-1)/2}{(p-3)/4}}\Big(1+p\Big(-1-\f 52\qp 2\Big)\Big)
\\&\e -6p^2(-1)^{\f{p+1}4}\b{(p-1)/2}{(p-3)/4}^{-1}(1+p\qp 2)^3
\Big(1+p\Big(-1-\f 52\qp 2\Big)\Big)
\\&\e
-6p^2(-1)^{\f{p+1}4}\b{(p-1)/2}{(p-3)/4}^{-1}\Big(1+p\Big(-1+\f
12\qp
2\Big)\Big)\\&=3p^2(-1)^{\f{p+1}4}\f{2p-1-2^{p-1}}{\b{(p-1)/2}{(p-3)/4}}
\mod {p^4}.\endalign$$
This completes the proof.

\pro{Theorem 5.4} Let $p>3$ be a prime. Then
$$\align &\sum_{k=0}^{p-1}(8k+3)\b{-3/4}k^4
\\&\e\cases (-1)^{\f{p-1}4}\big(\f{3p^2}{2x}+\f{3p^3}{8x^3}\big)
\mod {p^4}\q\t{if $p=x^2+y^2\e 1\mod 4$ and $4\mid x-1$,}
\\(-1)^{\f{p-3}4}p\b{\f{p-3}2}{\f{p-3}4}
\Big(1-p\big(1+\f 12\qp 2\big)+p^2\big(\f 12\qp 2+\f 38\qp 2^2+\f
14E_{p-3}\big)\Big) \mod {p^4}\\\qq\qq\qq\qq\qq\qq\qq\t{if $p\e
3\mod 4$.}
\endcases\endalign$$
\endpro
Proof. We first assume that $p\e 1\mod 4$ and so $p=x^2+y^2$ with
$4\mid x-1$. Set $a=-\f 34$. Then $\ap=\f{3p-3}4$ and $a'=-\f 34$.
By Theorem 5.1, Lemma 5.2 and (5.4),
$$\align Q_p\Big(-\f 34\Big)&\e -\f
34p^2(-16)^{\f{3(p-1)}4}\b{(p-1)/2}{(p-1)/4}^{-1} \Big(1+p\big(1-\f
52\qp 2+\f 14(H_{\f{p-1}4}-H_{\f{p-5}4})\Big)
\\&\e -\f 34p^2(-1)^{\f{p-1}4}\f{(1+p\qp 2)^3(1-\f 52p\qp 2)}{\f
12(2^{p-1}+1)(2x-\f p{2x})}
\\&\e -\f 34(-1)^{\f{p-1}4}\f{p^2}{2x-\f p{2x}}\e -\f 34(-1)^{\f{p-1}4}p^2
\f{2x+\f p{2x}}{4x^2} \mod {p^4}.
\endalign$$
Hence, $$\sum_{k=0}^{p-1}(8k+3)\b{-3/4}k^4=-4Q_p\Big(-\f 34\Big)\e
(-1)^{\f{p-1}4}\Big(\f{3p^2}{2x}+\f{3p^3}{8x^3}\Big) \mod {p^4}.$$
\par Now suppose that $p\e 3\mod 4$. For $a=-\f 34$ we have
$\ap=\f{p-3}4<\f p2$ and $a'=-\f 14$. By Theorem 5.1, Lemma 5.2 and
(5.4),
$$\align Q_p\Big(-\f 34\Big)&\e -\f p4\cdot
(-1)^{\f{p-3}4}\b{\f{p-3}2}{\f{p-3}4}\Big(1-\f
p2\big(H_{\f{p-3}2}-H_{\f{p-3}4}\big)\\&\q+\f{p^2}{16}\big(2(H_{\f{p-3}2}-H_{\f{p-3}4})^2
-2H_{\f{p-3}2}^{(2)}+H_{\f{p-3}4}^{(2)}\big)\Big)
\\&\e -\f p4\cdot
(-1)^{\f{p-3}4}\b{\f{p-3}2}{\f{p-3}4}\Big(1-\f p2\Big(\f 2{1-p}+\qp
2-\f p2\qp 2^2-pE_{p-3}\Big)\\&\q+\f{p^2}{16}\Big(2(2+\qp 2)^2+\f
2{(\f{p-1}2)^2}-4E_{p-3}\Big)\Big)
\\&\e -\f p4\cdot
(-1)^{\f{p-3}4}\b{\f{p-3}2}{\f{p-3}4} \Big(1-p\big(1+\f 12\qp
2\big)\\&\q+p^2\big(\f 12\qp 2+\f 38\qp 2^2+\f 14E_{p-3}\big)\Big)
\mod {p^4},
\endalign$$
which yields the result in the case $p\e 3\mod 4$. The proof is now
complete.

\pro{Lemma 5.3} Let $p>3$ be a prime. Then
$$\align&H_{[\f p3]}\e -\f 32\qp 3+\f 34p\qp 3^2-p\Ls p3U_{p-3}\mod
{p^2},
\\&H_{[\f p6]}\e -2\qp 2-\f 32\qp 3+p\Big(\qp 2^2+\f 34\qp
3^2\Big)-\f 52p\Ls p3U_{p-3}\mod {p^2},
\\&H_{[\f{2p}3]}\e -\f 32\qp 3+\f 34p\qp 3^2+2p\Ls p3U_{p-3}\mod
{p^2},
\\&H_{[\f p3]}^{(2)}\e -H_{[\f{2p}3]}^{(2)}\e 3\Ls p3U_{p-3}\mod p,
\\&H_{[\f p6]}^{(2)}\e 15\Ls p3U_{p-3}\mod p.
\endalign$$
\endpro
Proof. The first two congruences were given in [S4, Theorem 3.2]. By
[S4, Theorem 3.2], $\sum_{k=1}^{[2p/3]}\f{(-1)^{k-1}}k\e 3p\sls
p3U_{p-3}\mod{p^2}$. Thus,
$$H_{[\f{2p}3]}=H_{[\f p3]}+\sum_{k=1}^{[2p/3]}\f{(-1)^{k-1}}k
\e -\f 32\qp 3+\f 34p\qp 3^2+2p\Ls p3U_{p-3}\mod{p^2}.$$ By [S4,
Theorem 3.3], $H_{[\f p3]}^{(2)}\e 3\ls p3U_{p-3}\mod p$ and $H_{[\f
p6]}^{(2)}\e 15\sls p3U_{p-3}\mod p$. To complete the proof, we note
that
$$H_{[\f{2p}3]}^{(2)}=\sum_{k=1}^{p-1}\f 1{k^2}- \sum_{k=1}^{[p/3]}\f
1{(p-k)^2}\e -H_{[\f p3]}^{(2)}\mod p.$$

\pro{Theorem 5.5} Let $p>3$ be a prime. Then
$$\align &\sum_{k=0}^{p-1}(6k+1)\b{-1/3}k^4
\\&\e\cases -px+\f{p^2}x+\f{p^3}{x^3}\mod {p^4}
\q\t{if $3\mid p-1$ and so $4p=x^2+27y^2$ with $3\mid x-1$,}
\\p^2(-1)^{\f{p-1}2}(2^{p-1}-p)\cdot 2^{-\f{p-5}3}\b{\f{p-2}3}{\f{p-5}6}^{-1}
\mod {p^4}\q\t{if $p\e 2\mod 3$.}
\endcases\endalign$$
\endpro
Proof. We first assume that $p\e 1\mod 3$ and so $4p=x^2+27y^2$ with
$3\mid x-1$. Set $a=-\f 13$. Then $\ap=\f{p-1}3$ and $a'=-\f 13$. By
Theorem 5.1 and Lemma 5.3,
$$\align Q_p\Big(-\f 13\Big)&\e -\f
p3(-1)^{\f{p-1}3}\b{2(p-1)/3}{(p-1)/3}\Big(1-\f
23p\big(H_{\f{2(p-1)}3}-H_{\f{p-1}3}\big) \\&\q+\f{p^2}9
\big(2\big(H_{\f{2(p-1)}3}-H_{\f{p-1}3}\big)^2-2H_{\f{2(p-1)}3}^{(2)}+H_{\f{p-1}3}^{(2)}\big)
\Big)
\\&\e -\f
p3\b{2(p-1)/3}{(p-1)/3}\Big(1-\f 23p\cdot 3pU_{p-3}+\f {p^2}9\cdot
3\cdot 3U_{p-3}\Big)\\&=-\f
p3\b{2(p-1)/3}{(p-1)/3}\big(1-p^2U_{p-3}\big)\mod {p^4}.
\endalign$$
From [CD],
$$\b{2(p-1)/3}{(p-1)/3}\e \Big(-x+\f px+\f{p^2}{x^3}\Big)\Big(1+\f
16p^2B_{p-2}\Ls 13\Big)\mod {p^3}.$$ By [S4], $B_{p-2}\sls 13\e
6U_{p-3}\mod p$. Hence,
$$\align &\sum_{k=0}^{p-1}(6k+1)\b{-1/3}k^4\\&=-3Q_p\Big(-\f 13\Big)\e
p\Big(-x+\f px+\f{p^2}{x^3}\Big)\Big(1+\f 16p^2B_{p-2}\Ls
13\Big)\big(1-p^2U_{p-3}\big)
\\&\e p\Big(-x+\f px+\f{p^2}{x^3}\Big)\Big(1-(p^2U_{p-3})^2\Big)
\e p\Big(-x+\f px+\f{p^2}{x^3}\Big)\mod {p^4}.
\endalign$$
\par Now suppose that $p\e 2\mod 3$. For $a=-\f 13$ we have
$\ap=\f{2p-1}3>\f p2$ and $a'=-\f 23$. From Theorem 5.1 and Lemma
5.3 we deduce that
$$\align Q_p\Big(-\f 13\Big)&\e \f
23p^2(-1)^{\f{p-1}2}(-16)^{\f{2p-1}3}\b{\f{p+1}3}{\f{p+1}6}^{-1}
\Big(1+p\Big(-1-\f 73\qp 2+\f
16(H_{\f{p-2}3}-H_{\f{p-5}6})\Big)\Big)
\\&\e -\f
13p^2(-1)^{\f{p-1}2}2^{3(p-1)-\f{p-5}3}\b{(p-2)/3}{(p-5)/6}^{-1}
\\&\q\times\Big(1+p\Big(-1-\f 73\qp 2+\f
16\Big(-\f 32\qp 3+2\qp 2+\f 32\qp 3\Big)\Big)
\\&\e -\f 13p^2(-1)^{\f{p-1}2}(1+p\qp 2)^3
\cdot 2^{-\f{p-5}3}\b{(p-2)/3}{(p-5)/6}^{-1} (1-p(1+2\qp 2))
\\&\e -\f 13p^2(-1)^{\f{p-1}2}(2^{p-1}-p)
\cdot 2^{-\f{p-5}3}\b{(p-2)/3}{(p-5)/6}^{-1}\mod {p^4}.
\endalign$$
This yields the result in the case $p\e 2\mod 3$. The proof is now
complete.

\pro{Theorem 5.6} Let $p$ be a prime with $p>3$. If $p\e 1\mod 3$
and so $p=x^2+3y^2$ with $x\e 1\mod 3$ and $y\e 0,1\mod 3$, setting
$$u=\cases -2x&\t{if $3\mid y$,}
\\x-3y&\t{if $3\mid y-1$}
\endcases\tag 5.5$$
we then have
$$\sum_{k=0}^{p-1}(12k+1)\b{-1/6}k^4\e -pu+\f{p^2}u+\f{p^3}{u^3}\mod
{p^4}.$$ If $p\e 2\mod 3$, then
$$\sum_{k=0}^{p-1}(12k+1)\b{-1/6}k^4\e
\f{20\cdot
2^{\f{p+1}3}(2^{p-1}+3p-4)p^2}{3\b{2(p+1)/3}{(p+1)/3}}\mod {p^4}.$$
\endpro
Proof. We first assume that $p\e 1\mod 3$ and so $p=x^2+3y^2$ with
$x\e 1\mod 3$ and $y\e 0,1\mod 3$. For $a=-\f 16$ we have
$\ap=\f{p-1}6$ and $a'=-\f 16$. By Theorem 5.1 and Lemma 5.3,
$$\align Q_p\Big(-\f 16\Big)&\e -\f
p6(-1)^{\f{p-1}6}\b{(p-1)/3}{(p-1)/6}\Big(1-\f
p3\big(H_{\f{p-1}3}-H_{\f{p-1}6}\big)
\\&\q+\f{p^2}{36}\Big(2\big(H_{\f{p-1}3}-H_{\f{p-1}6}\big)^2-2H_{\f{p-1}3}^{(2)}
+H_{\f{p-1}6}^{(2)}\Big)\Big)
\\&\e -\f p6(-1)^{\f{p-1}6}\b{(p-1)/3}{(p-1)/6}\Big(1-\f
p3\big(2\qp 2-p\qp 2^2+\f 32pU_{p-3}\big)
\\&\q+\f{p^2}{36}\big(2(2\qp 2)^2-6U_{p-3}+15U_{p-3}\big)\Big)
\\&=-\f p6(-1)^{\f{p-1}6}\b{(p-1)/3}{(p-1)/6}\Big(1-\f 23p\qp
2+p^2\big(\f 59\qp 2^2-\f 14U_{p-3}\big)\Big)\mod {p^4}.
\endalign$$
By [SD],
$$\b{\f{p-1}3}{\f{p-1}6}\e (-1)^{\f{p-1}6-1}\Big(u-\f
pu-\f{p^2}{u^3}\Big)\Big(1+\f 23p\qp 2+p^2\big(-\f 19\qp 2^2+\f
1{24}B_{p-2}\ls 13\big)\Big)\mod {p^3},$$ where $u$ is given by
(5.5). Since $B_{p-1}\sls 13\e 6U_{p-3}\mod p$ by [S4], from the
above we deduce that
$$\align -6Q_p\Big(-\f 16\Big)&\e -p\Big(u-\f pu-\f{p^2}{u^3}\Big)\Big(1+\f
23p\qp 2+p^2\big(-\f 19\qp 2^2+\f 14U_{p-3}\big)\Big)
\\&\q\times \Big(1-\f 23p\qp
2+p^2\big(\f 59\qp 2^2-\f 14U_{p-3}\big)\Big)
\\&\e -p\Big(u-\f pu-\f{p^2}{u^3}\Big)\mod {p^4}.
\endalign$$
This yields the result in the case $p\e 1\mod 3$.
\par Now assume that $p\e 2\mod 3$. For $a=-\f 16$ we see that
$\ap=\f{5p-1}6>\f p2$ and $a'=-\f 56$. By Lemma 5.3,
$$H_{p-1-\ap}-H_{\ap-\f{p+1}2}=H_{\f{p-5}6}-H_{\f{p-2}3}
\e -2\qp 2-\f 32\qp 3+\f 32\qp 3=-2\qp 2\mod p.$$ Hence, taking
$a=-\f 16$ in Theorem 5.1 gives
$$\align Q_p\Big(-\f 16\Big)&\e
\f 56p^2(-1)^{\f{p-1}2}(-16)^{\f{5p-1}6}\b{2(p+1)/3}{(p+1)/3}^{-1}
\Big(1+p\big(-1-\f 83\qp 2-\f 23\qp 2\big)\Big)
\\&\e \f 56 p^2(1+p\qp 2)^3\cdot 2^{\f{p+7}3}\b{2(p+1)/3}{(p+1)/3}^{-1}
\Big(1-p\big(1+\f {10}3\qp 2\big)\Big)
\\&\e \f 56p^2\cdot 4\cdot
2^{\f{p+1}3}\b{2(p+1)/3}{(p+1)/3}^{-1}\Big(1-p\big(1+\f 13\qp
2\big)\Big)\mod {p^4}.
\endalign$$
This yields the result in the case $p\e 2\mod 3$. The proof is now
complete.

\pro{Theorem 5.7} Let $p>3$ be a prime. If $p\e 1,3\mod 8$ and so
$p=x^2+2y^2$, then
$$\align&\sum_{k=0}^{p-1}(16k+1)\b{-1/8}k^4
\\&\e\cases (-1)^{\f
y2}p\big(2x-\f{p}{2x}-\f{p^2}{8x^3}\big)\mod {p^4}&\t{if $8\mid p-1$
and $4\mid x-1$,}\\3p\big(4y-\f{p}{2y}-\f{p^2}{16y^3}\big)\mod
{p^4}&\t{if $8\mid p-3$ and $4\mid y-1$.}
\endcases\endalign$$
If $p\e 5,7\mod 8$, then
$$\sum_{k=0}^{p-1}(16k+1)\b{-1/8}k^4
\e\cases (-1)^{\f{p+3}8}\f{20p^2}{3\binom{(p+3)/4}{(p+3)/8}} \mod
{p^3}&\t{if $8\mid p-5$,}
\\(-1)^{\f{p-7}8}\f{56p^2}{\b{3(p+1)/4}{3(p+1)/8}}
\mod {p^3}&\t{if $8\mid p-7$.}\endcases$$
\endpro
Proof. We first assume that $p\e 1\mod 8$ and so $p=x^2+2y^2$ with
$4\mid x-1$. For $a=-\f 18$ we have $\ap=\f{p-1}8<\f p2$ and $a'=-\f
18$. Taking $a=-\f 18$ in Theorem 5.1 gives
$$\f{Q_p(-\f 18)}p\e -\f 18(-1)^{\f{p-1}8}\b{(p-1)/4}{(p-1)/8}
 \mod {p}.$$
Clearly, $\b{(p-1)/2}{(p-1)/8}\e
\b{-1/2}{(p-1)/8}=\b{(p-1)/4}{(p-1)/8}(-4)^{-\f{p-1}8}\mod p$. By a
result due to Stern (see [BEW]),
$$\b{(p-1)/2}{(p-1)/8}\e (-1)^{\f{p-1}8}\cdot
2x\mod p.$$ Thus,
$$ \b{(p-1)/4}{(p-1)/8}\e 4^{\f{p-1}8}\cdot
2x=2^{\f{p-1}4}\cdot 2x\e (x/y)^{\f{p-1}2} \e\Ls {xy}p\cdot 2x\mod
p.$$ Suppose $y=2^{\alpha}y_0\;(2\nmid y_0)$. Then $\sls{y}p=\sls
2p^{\alpha}\sls {y_0}p=\sls{y_0}p=\sls p{|y_0|}=\sls{x^2}{|y_0|}=1$
and $\sls xp=\sls p{|x|}=\sls{x^2+2y^2}{|x|}=\sls
2{|x|}=(-1)^{\f{x^2-1}8}=(-1)^{\f{p-1-2y^2}8}=(-1)^{\f{p-1}8+\f
y2}$. Hence,
$$\f{-8Q_p(-\f 18)}p\e (-1)^{\f{p-1}8}\b{\f{p-1}4}{\f{p-1}8}
\e (-1)^{\f{p-1}8}\Ls xp\Ls yp\cdot 2x=(-1)^{\f y2}\cdot 2x\mod p.$$
From Theorem 5.2 we have
$$\Ls{Q_p(-\f 18)}p^2\e 2^{\f{p-1}2}\cdot\f 1{64}S_p\Big(-\f 14\Big)
\Big(1-\f 12p\qp 2+\f 38p^2\qp 2^2\Big)\mod {p^3}.$$
 By [M],
 $$S_p\Big(-\f
 14\Big)=\sum_{k=0}^{p-1}\f{\b{2k}k^2\b{4k}{2k}}{256^k}
 \e 4x^2-2p-\f{p^2}{4x^2}\mod {p^3}.$$
 It is easy to see that
 $$\Ls 2p2^{\f{p-1}2}\e 1+\f 12p\qp 2-\f 18p^2\qp 2^2\mod {p^3}.\tag 5.6$$
 Thus,
 $$\align\Ls{-8Q_p(-\f 18)}p^2&\e \Big(4x^2-2p-\f{p^2}{4x^2}\Big)
 \Big(1+\f 12p\qp 2-\f 18p^2\qp 2^2\Big)
\\&\q\times\Big(1-\f 12p\qp 2+\f
38p^2\qp 2^2\Big)
\\&\e 4x^2-2p-\f{p^2}{4x^2}
\e \Big(2x-\f{p}{2x}-\f{p^2}{8x^3}\Big)^2\mod {p^3}.
\endalign$$
That is, $$\Big(\f{-8Q_p(-\f 18)}p-(-1)^{\f
y2}\Big(2x-\f{p}{2x}-\f{p^2}{8x^3}\Big)\Big) \Big(\f{-8Q_p(-\f
18)}p+(-1)^{\f y2}\Big(2x-\f{p}{2x}-\f{p^2}{8x^3}\Big)\Big) \e 0\mod
{p^3}.$$ Since $\f{-8Q_p(-\f 18)}p\e (-1)^{\f
y2}(2x-\f{p}{2x}-\f{p^2}{8x^3})\mod p$, we must have
$$\f{-8Q_p(-\f 18)}p\e(-1)^{\f
y2}\Big(2x-\f{p}{2x}-\f{p^2}{8x^3}\Big)\mod {p^3}.$$ This yields the
result in the case $p\e 1\mod 8$.
\par Now suppose that $p\e 3\mod
8$ and so $p=x^2+2y^2$ with $x\e y\e 1\mod 4$. For $a=-\f 18$ we
have $\ap=\f{3p-1}8<\f p2$ and $a'=-\f 38$. Taking $a=-\f 18$ in
Theorem 2.1 gives
$$\align \f{Q_p(-\f 18)}p&\e -\f 38(-1)^{\f{3p-1}8}\b{(3p-1)/4}{(3p-1)/8}
=-\f 38\cdot 4^{\f{3p-1}8}\b{-1/2}{(3p-1)/8}
 \\&\e -\f 38\cdot 2^{\f{3p-1}4}\b{(p-1)/2}{(3p-1)/8}
\e -\f 38\b{(p-1)/2}{(p-3)/8}(-2)^{\f{p+1}4}
\\&\e -\f 38\b{(p-1)/2}{(p-3)/8}\ls xy^{\f{p+1}2}
\e -\f 38\b{(p-1)/2}{(p-3)/8}\f xy\Ls xp\Ls yp
 \mod {p}.\endalign$$
 Observe that
 $$\align&\Ls xp=\Ls p{|x|}=\Ls{x^2+2y^2}{|x|}=\Ls 2{|x|}
 =(-1)^{\f{x^2-1}8}=(-1)^{\f{p-3-2(y^2-1)}8}=(-1)^{\f{p-3}8},
\\&\Ls yp=\Ls p{|y|}=\Ls{x^2+2y^2}{|y|}=\Ls{x^2}{|y|}=1.\endalign$$
By a result due to Eisenstein (see [BEW]),
$$\b{(p-1)/2}{(p-3)/8}\e -2(-1)^{\f{p-3}8}x\mod p.$$
Hence,
$$\f{-8Q_p(-\f 18)}p\e 3\b{(p-1)/2}{(p-3)/8}\f xy\Ls xp\Ls yp
\e 3(-1)^{\f{p-3}8}(-2x)\cdot\f xy\cdot (-1)^{\f{p-3}8}\e 12y\mod
p.$$ Taking $a=-\f 18$ in Theorem 5.2 yields
$$\Ls{Q_p(-\f 18)}p^2\e 2^{\f{3p-1}2}\Big(-\f 38\Big)^2S_p\Big(-\f
14\Big)\Big(1-\f 32p\qp 2+\f{15}8p^2\qp 2^2\Big)\mod {p^3}.$$ By
(5.6),
$$\align 2^{\f{3(p-1)}2}&=2^{p-1}\cdot 2^{\f{p-1}2}
\e -(1+p\qp 2)\Big(1+\f 12p\qp 2-\f 18p^2\qp 2^2\Big)
\\&\e -\Big(1+\f 32p\qp
2+\f 38p^2\qp 2^2\Big)\mod {p^3}.\endalign$$ By [M], $S_p(-\f 14)\e
4x^2-2p-\f{p^2}{4x^2}\mod {p^3}$. Hence,
$$\align\Ls{-8Q_p(-\f 18)}p^2&\e 9\Big(4x^2-2p-\f{p^2}{4x^2}\Big)
\cdot (-2)\Big(1+\f 32p\qp 2+\f 38p^2\qp 2^2\Big)
\\&\q\times \Big(1-\f 32p\qp 2+\f{15}8p^2\qp 2^2\Big)
\\&\e -18\Big(4x^2-2p-\f{p^2}{4x^2}\Big)\mod {p^3}.\endalign$$
It is clear that
$$\align \Big(4y-\f p{2y}-\f{p^2}{16y^3}\Big)^2
&\e 16y^2-4p-\f{p^2}{4y^2}=8(p-x^2)-4p-\f{p^2}{2(p-x^2)}
\\&\e -2\Big(4x^2-2p-\f{p^2}{4x^2}\Big)\mod {p^3}.\endalign$$
Thus,
$$\align &\Big(\f{-8Q_p(-\f 18)}p-3\Big(4y-\f
p{2y}-\f{p^2}{16y^3}\Big)\Big)\Big(\f{-8Q_p(-\f 18)}p+3\Big(4y-\f
p{2y}-\f{p^2}{16y^3}\Big)\Big)
\\&\e\Ls{-8Q_p(-\f 18)}p^2+18\Big(4x^2-2p-\f{p^2}{4x^2}\Big)\e 0\mod
{p^3}.\endalign$$ Since $$\f{-8Q_p(-\f 18)}p\e 12y\e 3\Big(4y-\f
p{2y}-\f{p^2}{16y^3}\Big)\mod p,$$ we must have
$$\f{-8Q_p(-\f 18)}p\e  3\Big(4y-\f
p{2y}-\f{p^2}{16y^3}\Big)\mod {p^3}.$$ This proves the result in the
case $p\e 3\mod 8$.
\par For $p\e 5\mod 8$ taking $\ap=\f{5p-1}8$ and $a'=-\f 58$ in
Theorem 5.1 yields the result. For $p\e 7\mod 8$ taking
$\ap=\f{7p-1}8$ and $a'=-\f 78$ in Theorem 5.1 yields the result.
\par Summarizing the above proves the theorem.
\par\q
\par In a similar way, one can prove the following result.
 \pro{Theorem 5.8} Let $p>3$ be a prime. If $p\e 1,3\mod 8$ and so
$p=x^2+2y^2$, then
$$\align&\sum_{k=0}^{p-1}(16k+3)\b{-3/8}k^4
\\&\e\cases 3p(-1)^{\f
y2}\big(2x-\f{p}{2x}-\f{p^2}{8x^3}\big)\mod {p^4}&\t{if $8\mid p-1$
and $4\mid x-1$,}\\p\big(-2y+\f{p}{4y}+\f{p^2}{32y^3}\big)\mod
{p^4}&\t{if $8\mid p-3$ and $4\mid y-1$.}
\endcases\endalign$$
If $p\e 5,7\mod 8$, then
$$\sum_{k=0}^{p-1}(16k+3)\b{-3/8}k^4
\e\cases (-1)^{\f{p-5}8}\f{84p^2}{\binom{(3p+1)/4}{(3p+1)/8}} \mod
{p^3}&\t{if $8\mid p-5$,}
\\(-1)^{\f{p-7}8}\f{10p^2}{\b{(p+1)/4}{(p+1)/8}}
\mod {p^3}&\t{if $8\mid p-7$.}\endcases$$
\endpro

\pro{Theorem 5.9} Let $p>3$ be a prime. If $p\e 1\mod 4$ and so
$p=x^2+y^2$ with $4\mid x-1$, then
$$\align&\sum_{k=0}^{p-1}(24k+1)\b{-1/12}k^4
\\&\e\cases 2^{\f{p-1}6}p\big(1-\f 16(2^{p-1}-1)\big)\big(2x-\f p{2x}\big)\mod {p^3}
&\t{if $12\mid p-1$ and $3\mid
x$,}
\\-2^{\f{p-1}6}p\big(1-\f 16(2^{p-1}-1)\big)\big(2x-\f p{2x}\big)\mod {p^3}
&\t{if $12\mid p-1$ and $3\nmid x$,}
\\5\cdot 2^{-\f{p-5}6}p\big(1+\f 16(2^{p-1}-1)\big)\big(2y-\f p{2y}\big)\mod {p^3}&\t{if $12\mid p-5$ and $3\mid x-y$.}
\endcases\endalign$$
If $p\e 3\mod 4$, then
$$\align\sum_{k=0}^{p-1}(24k+1)\b{-1/12}k^4
\e\cases (-1)^{\f{p-7}{12}}\f{28\cdot 2^{\f{p-1}3}p^2}
{5\b{(p+5)/6}{(p+5)/12}}\mod {p^3}&\t{if $12\mid p-7$,}
\\(-1)^{\f{p-11}{12}}\f{11p^2}{2^{\f{p-11}3}\b{5(p+1)/6}{5(p+1)/12}}\mod {p^3}&\t{if $12\mid p-11$.}
\endcases\endalign$$

\endpro
Proof. We first assume that $p\e 1\mod 4$ and so $p=x^2+y^2$ with
$4\mid x-1$. By [S7, p.11],
$$\b{(p-1)/2}{[p/12]}\e \Big(2c-\f p{2c}\Big)\Big(1+p\Big(\f 32q_p(2)
+\f 54q_p(3)+\f 13H_{[\f p{12}]}\Big)\Big)\mod {p^2},\tag 5.7$$
where
$$c=\cases x&\t{if $p\e 1\mod {12}$ and $3\nmid x$,}
\\-x&\t{if $p\e 1\mod {12}$ and $3\mid x$,}
\\y&\t{if $p\e 5\mod {12}$ and $y\e x\mod 3$.}
\endcases$$
From [S3, Lemma 2.4], for $k=1,2,\ldots,\f{p-1}2$,
$$\b{2k}k\e \b{(p-1)/2}k(-4)^k\Big(1+p\big(H_{2k}-\f 12H_k\big)\Big)
\mod {p^2}.\tag 5.8$$

For $p\e 1\mod {12}$ and $a=-\f 1{12}$ we have $\ap=\f{p-1}{12}$ and
$a'=-\f 1{12}$. From the above, Lemma 5.3 and Theorem 5.1 we deduce
that
$$\align\f{-12Q_p(-\f 1{12})}p&\e (-1)^{\f{p-1}{12}}\b{(p-1)/6}{(p-1)/12}\Big(1-\f
p6\big(H_{\f{p-1}6}-H_{\f{p-1}{12}}\big)\Big)
\\&\e
4^{\f{p-1}{12}}\b{(p-1)/2}{(p-1)/12}\big(1+p\big(H_{\f{p-1}6}-\f
12H_{\f{p-1}{12}}\big)\big)\Big(1-\f
p6\big(H_{\f{p-1}6}-H_{\f{p-1}{12}}\big)\Big)
\\&\e 2^{\f{p-1}6}\b{(p-1)/2}{(p-1)/12}\Big(1+p\Big(\f
56H_{\f{p-1}6}-\f 13H_{\f{p-1}{12}}\Big)\Big)
\\&\e 2^{\f{p-1}6}\b{(p-1)/2}{(p-1)/12}\Big(1+p\Big(-\f
53\qp 2-\f 54\qp 3-\f 13H_{\f{p-1}{12}}\Big)\Big) \\&\e
2^{\f{p-1}6}\Big(2c-\f p{2c}\Big)\Big(1+p\Big(\f 32q_p(2) +\f
54q_p(3)+\f 13H_{[\f p{12}]}\Big)\Big)\\&\q\times \Big(1+p\Big(-\f
53\qp 2-\f 54\qp 3-\f 13H_{\f{p-1}{12}}\Big)\Big)
\\&\e 2^{\f{p-1}6}\big(1-\f 16(2^{p-1}-1)\big)\big(2c-\f p{2c}\big)\mod
{p^3}.\endalign$$ For $p\e 5\mod {12}$ and $a=-\f 1{12}$ we have
$\ap=\f{5p-1}{12}$ and $a'=-\f 5{12}$. From (5.7), (5.8) and Theorem
5.1 we deduce that
$$\align\f{-12Q_p(-\f 1{12})}p&\e 5(-1)^{\f{5p-1}{12}}\b{(5p-1)/6}{(5p-1)/12}
\Big(1-\f 56p\big(H_{\f{5p-1}6}-H_{\f{5p-1}{12}}\big)\Big)
\\&\e 5\cdot 4^{\f{5p-1}{12}}\b{(p-1)/2}{(5p-1)/12}
\Big(1+p\big(H_{\f{5p-1}6}-\f 12H_{\f{5p-1}{12}}\big)\Big) \Big(1-\f
56p\big(H_{\f{5p-1}6}-H_{\f{5p-1}{12}}\big)\Big)
\\&\e 5\cdot 2^{p-1-\f{p-5}6}\b{(p-1)/2}{(p-5)/12}\Big(1+p\Big(\f
16H_{\f{5p-1}6}+\f 13H_{\f{5p-1}{12}}\Big)\Big)\mod {p^3}.
\endalign$$
By [S7, (2.1) and (2.3)], for $k\in\{1,2,\ldots,\f{p-1}2\}$ we have
$$H_{p-1-k}\e H_k\mod p\qtq{and} H_k+H_{\f{p-1}2-k}\e 2H_{2k}-2\qp 2\mod p.$$
Hence, appealing to Lemma 5.3 we get
$$\align \f 16H_{\f{5p-1}6}+\f 13H_{\f{5p-1}{12}}&\e \f
16H_{\f{p-5}6}+\f 13\Big(2H_{\f{p-5}6}-2\qp 2-H_{\f{p-5}{12}}\Big)
\\&\e \f 56\Big(-2\qp 2-\f 32\qp 3\Big)-\f 23\qp 2-\f
13H_{\f{p-5}{12}}
\\&=-\f 73\qp 2-\f 54\qp 3-\f 13H_{\f{p-5}{12}}\mod {p^3}.
\endalign$$
Therefore,
$$\align \f{-12Q_p(-\f 1{12})}p&\e 5\cdot 2^{-\f{p-5}6}(1+p\qp
2)\b{\f{p-1}2}{\f{p-5}{12}}\Big(1+p\Big(-\f 73\qp 2-\f 54\qp 3-\f
13H_{\f{p-5}{12}}\Big)\Big)
\\&\e 5\cdot 2^{-\f{p-5}6}(1+p\qp 2)\Big(2c-\f p{2c}\Big)\Big(1+p\Big(\f 32\qp 2+\f 54\qp 3+\f 13H_{\f{p-5}{12}}\Big)\Big)
\\&\q\times \Big(1+p\Big(-\f 73\qp 2-\f 54\qp 3-\f
13H_{\f{p-5}{12}}\Big)\Big)
\\&\e 5\cdot 2^{-\f{p-5}6}\Big(1+\f 16p\qp 2\Big)\Big(2c-\f p{2c}\Big)\mod {p^3}.
\endalign$$
\par For $p\e 7\mod {12}$ and $a=-\f 1{12}$ we have $\ap=\f{7p-1}{12}$
and $a'=-\f 7{12}$. By Theorem 5.1,
$$ -12Q_p(-\f 1{12})\e \f
75p^2(-16)^{\f{7p-1}{12}}\b{\f{p+5}6}{\f{p+5}{12}}^{-1}\e \f
{28}5p^2(-1)^{\f{p-7}{12}}2^{\f{p-1}3}\b{\f{p+5}6}{\f{p+5}{12}}^{-1}
\mod {p^3}.$$ For $p\e 11\mod {12}$ and $a=-\f 1{12}$ we have
$\ap=\f{11p-1}{12}$ and $a'=-\f {11}{12}$. By Theorem 5.1,
$$ -12Q_p(-\f 1{12})\e 11p^2(-16)^{\f{11p-1}{12}}\b{\f{5(p+1)}6}
{\f{5(p+1)}{12}}^{-1}\e
11p^2(-1)^{\f{p-11}{12}}2^{-\f{p-11}3}\b{\f{5(p+1)}6}
{\f{5(p+1)}{12}}^{-1}\mod {p^3}.$$ This completes the proof.

\end{document}